\documentclass[11pt]{amsart}

\usepackage{amsmath,amssymb,amsthm,mathtools}
\usepackage{booktabs,array}
\usepackage{mathrsfs}
\usepackage[expansion=false]{microtype}
\usepackage[hidelinks]{hyperref}
\usepackage[margin=1.25in]{geometry}

\newtheorem{theorem}{Theorem}[section]
\newtheorem{proposition}[theorem]{Proposition}

\newtheorem{corollary}[theorem]{Corollary}
\theoremstyle{definition}

\theoremstyle{remark}
\newtheorem{remark}[theorem]{Remark}

\newcommand{\Q}{\mathbf Q}
\newcommand{\Z}{\mathbf Z}
\newcommand{\C}{\mathbf C}
\newcommand{\F}{\mathbf F}

\newcommand{\Jac}{\operatorname{Jac}}
\newcommand{\Gal}{\operatorname{Gal}}
\newcommand{\Tr}{\operatorname{Tr}}
\newcommand{\im}{\operatorname{im}}

\newcommand{\PSL}{\operatorname{PSL}}
\newcommand{\PGL}{\operatorname{PGL}}
\newcommand{\SL}{\operatorname{SL}}
\newcommand{\GL}{\operatorname{GL}}
\newcommand{\diag}{\operatorname{diag}}

\title[Hypergeometric abelian surfaces]{Hypergeometric Abelian Surfaces and Triangle-Group Levels}
\author{Devendra Tiwari}
\address{Statistics and Mathematics Department, Indian Statistical Institute, Bangalore 560 059, INDIA}
\email{devendra9.dev@gmail.com}
\subjclass[2020]{14K02, 14H40, 20H10, 33C05}
\keywords{hypergeometric abelian varieties, Frey abelian varieties, cyclic covers, real multiplication, triangle groups, triangular modular curves}
\date{November 2026}

\begin{document}

\begin{abstract}
We study three families of hypergeometric abelian surfaces with real
multiplication by $K=\Q(\sqrt5)$. We determine their integral homology
lattices and monodromy, including the choice of lattice at the ramified
prime above $5$ in the plus family. One family is isomorphic to an explicit
quadratic twist of a genus-two Frey Jacobian, with the two real-multiplication
embeddings exchanged. The specified isogeny is twice an isomorphism and
has kernel equal to the full $2$-torsion. Each family gives two real
rank-two representations, with projective local monodromies of order $5$
at two punctures and parabolic monodromy at the third. Exactly one of
the six projective representations is discrete, with image
$\Delta(5,5,\infty)$. Reduction at the two primes above $11$ gives six
degree-$12$, genus-two covers parametrizing torsion lines. For the family
with a discrete constituent, these Borel covers are the triangular modular
curves $X_0(5,5,\infty;\mathfrak q)$; the comparison descends to
ordered covers over $K$. The generic arithmetic projective image is
$\PGL_2(\F_{11})$; the Galois closure of each Borel cover has constant
field $K(\sqrt{-11})$.
\end{abstract}

\maketitle

\section{Introduction and main results}

\subsection{The geometric problem}

The Legendre family illustrates a familiar passage from geometry to level structure.  An algebraic family carries a variation of cohomology, its periods satisfy a hypergeometric differential equation, and its integral homology makes it possible to impose finite level conditions.  Hypergeometric abelian varieties provide a higher-dimensional setting for the same circle of ideas: their periods are carried by abelian factors of cyclic-cover Jacobians, while the finite covers attached to them depend not only on rational monodromy but on the integral lattice supplied by the geometry.

We recall the triangle-group notation that organizes this picture.  For
\[
 a,b,c\in\{2,3,\ldots\}\cup\{\infty\},
\]
the oriented triangle orbifold group is
\[
 \Delta(a,b,c)
 =\left\langle \gamma_0,\gamma_1,\gamma_\infty
 \,\middle|\,
 \gamma_0^a=\gamma_1^b=\gamma_\infty^c=
 \gamma_0\gamma_1\gamma_\infty=1
 \right\rangle,
\]
with a power relation omitted whenever its corresponding order is
$\infty$. Geometrically this is the orbifold fundamental group of a sphere
with three distinguished points; an order $\infty$ corresponds to a cusp.
A Gauss hypergeometric equation has singular points $0,1,\infty$.
A nonintegral rational local exponent difference with denominator $a$ in
lowest terms gives projective local monodromy of order $a$. In the systems
considered here, exponent difference $0$ gives nontrivial parabolic
monodromy, recorded by $\infty$. Thus the hypergeometric local system gives
a projective representation of the abstract group $\Delta(a,b,c)$ determined
by its local orders. This representation need not be faithful or discrete;
determining when its image is a Fuchsian triangle group is one of our
questions.

Darmon's construction associates hypergeometric abelian varieties to
generalized Fermat equations. For signature $(r,r,p)$, the complex
rank-two local systems have projective local monodromies of order $r$
at two punctures and parabolic monodromy at the third
\cite{Darmon2000,DarmonLecture}. The exponent $p$ is the characteristic
of the residual representation used in the arithmetic specialization,
whereas the local order $\infty$ describes the characteristic-zero
monodromy. The geometric dictionary in the Fourteenth Lecture suggests
studying these constructions together. For each abelian family we identify
the relevant rank-two constituents, recover their integral homology, and
compare the resulting finite covers with triangle-group congruence covers
when a discrete real image is available.

We work in the first abelian-surface case of the $(r,r,\infty)$ branch, namely $r=5$.  The three surfaces are
\[
 J^-_{5,5}(1),\qquad J^-_{5,5}(3),\qquad J^+_{5,5}(1),
\]
all with real multiplication by
\[
 K=\Q(\sqrt5),\qquad R=\mathcal O_K=\Z[\alpha],\qquad
 \alpha=\frac{\sqrt5-1}{2},\qquad \phi=1+\alpha=\frac{1+\sqrt5}{2}.
\]
These families arise from Darmon's cyclic-cover construction. We abbreviate
real multiplication by RM. For the minus sign, the index $j$ ranges
through $1,3,\ldots,r-2$; we call the last value the \emph{endpoint}.
Thus the endpoint family in this paper is $J^-_{5,5}(3)$.
The three families need not have the same global monodromy. We follow
each rank-two $K$-local system through its algebraic realization,
integral homology, two real embeddings, and finite level covers.

Throughout, $U=\mathbf P^1\setminus\{0,1,\infty\}$ is the open base of
the families, with the field of definition indicated by a subscript when needed.
A \emph{marked cover} is a branched cover of the fixed projective line
with branch values ordered as $(0,1,\infty)$; its restriction to $U$
is unramified, and its sheets are not individually labelled.
We call a cover \emph{regular} if it is Galois over an algebraic closure
of the constant field. A \emph{deck-labelled regular cover} additionally
has a specified identification of its geometric deck group with the
stated finite group.
Geometric monodromy is computed over an algebraic closure of the constant
field, whereas arithmetic monodromy also includes its Galois action.

The RM action makes integral first homology an $R$-lattice of rank two.
These lattices vary as a local system over $U$. A lattice is
\emph{monodromy-stable} if continuation along every loop preserves it.
A $K$-rational matrix representation alone does not identify the lattice
coming from the abelian variety. In the plus family, two homothety
classes of local monodromy-stable lattices occur at the prime above $5$;
the geometric construction selects one. Here homothety identifies
lattices that differ by multiplication by a nonzero element of the
coefficient field.

\subsection{Main results}

Our first result relates $J^-_{5,5}(3)$ to the genus-two Frey Jacobian
used in the corresponding generalized Fermat problem.

\begin{theorem}[Endpoint hyperelliptic comparison]\label{thm:intro-kraus}
Let $C'_5(t)$ denote the family of genus-two curves over $U$ described in Section~\ref{sec:r5kraus}. Then
\[
 J^-_{5,5}(3)\sim_{K(t)}\Jac(C'_5(t))^{[b(1-t)]},
 \qquad b=-\frac{5+\sqrt5}{2}.
\]
If $i_D$ and $i_B$ denote the real-multiplication actions on the two sides, the isogeny $\Psi$ satisfies
\[
 \Psi\,i_D(a)=i_B(\bar a)\,\Psi\qquad(a\in K),
\]
where $a\mapsto\bar a$ is the nontrivial automorphism of $K$.
The comparison therefore specifies both the quadratic twist and the change
of RM embedding.
\end{theorem}

The second group of results concerns the integral and real monodromy.  The lattice statement is stronger than the observation that one can write integral matrices: it identifies the homology lattice inherited from the abelian variety inside the fixed rational monodromy space.

\begin{theorem}[Integral homology lattices]\label{thm:intro-real-lattice}
For each of
\[
 J^-_{5,5}(1),\qquad J^-_{5,5}(3),\qquad J^+_{5,5}(1),
\]
the $R$-lattice inherited from integral first homology is $R^2$, up to
$K^\times$-homothety, in the $K$-basis of Section~\ref{sec:monodromy}.
In the plus family the integral action of the quotient involution selects
this lattice from the two local homothety classes at the prime above $5$.

For the plus Jacobian, after this identification, the canonical principal pairing has the form
\[
 E(v,w)=\Tr_{K/\Q}\left(\frac{\varepsilon}{\sqrt5}\det_K(v,w)\right),
 \qquad \varepsilon\in R^\times.
\]
\end{theorem}

Applying the two real embeddings of $K$ to the monodromy of the three surfaces produces six real projective representations.  Their behaviour is sharply nonuniform.

\begin{proposition}[Real projective images]\label{thm:intro-real-images}
Let $\sigma_+(\sqrt5)>0$ and $\sigma_-(\sqrt5)<0$.  Exactly one of these six representations is discrete: the representation attached to $J^-_{5,5}(3)$ under $\sigma_-$ is faithful with image
\[
 \Delta(5,5,\infty)\cong C_5*C_5.
\]
The other five are nondiscrete.
\end{proposition}

We next reduce these integral representations modulo the primes above
$11$. Choosing a line in the two-dimensional RM torsion component gives
a degree-$12$ cover. A line stabilizer is a Borel subgroup, and we call
this the \emph{Borel cover}.

\begin{theorem}[Level $11$]\label{thm:intro-level11}
For each of the three abelian surfaces and each of the two primes $\mathfrak p\mid11$ of $K$, the geometric monodromy reduces to $\SL_2(\F_{11})$, with projective image $\PSL_2(\F_{11})$.  The associated Borel cover has degree $12$ and genus $2$.  With the branch values $0,1,\infty$ ordered, the six resulting covers are pairwise inequivalent as marked covers.

For $J^-_{5,5}(3)$, whose $\sigma_-$ constituent is discrete, the two Borel covers are the triangular modular curves
\[
 X_0(5,5,\infty;\mathfrak q)
\]
for the two primes $\mathfrak q\mid11$ in the trace field of
$\Delta(5,5,\infty)$. The field identification in
Proposition~\ref{prop:CV11-endpoint} exchanges the two prime labels.
The four covers attached to the other two surfaces are not these two
Clark--Voight covers under the marking changes specified below.
\end{theorem}

These results relate the cyclic-cover construction, the Frey Jacobian,
integral RM monodromy, and finite level covers. The local projective
orders $(5,5,\infty)$ do not determine a discrete global image: only one
of the six real representations is discrete, although every reduction
at level $11$ has projective image $\PSL_2(\F_{11})$.

Corollary~\ref{cor:reverse20} refines the endpoint comparison:
the constructed map is twice a $K(t)$-isomorphism, has degree $16$
and kernel $J^-_{5,5}(3)[2]$, and admits a reverse morphism
whose compositions are multiplication by $20$.  Its integral
Tate-module comparison therefore holds away from $2$, including
at the coefficient prime above $5$.
We also determine the generic arithmetic image at level $11$ and the fields
of definition of the marked Borel and regular covers.  These arithmetic
statements concern the families over $K(t)$; they do not assert full
Galois image for every specialization or a Diophantine lifting theorem.

\subsection{Level structures and triangular modular curves}

The degree-$12$ covers in Theorem~\ref{thm:intro-level11} have a direct interpretation relative to the abelian family.  Let $\mathcal A\to U$ be one of the smooth RM families and let $\Lambda$ denote its rank-two integral homology local system over $R$. For either coefficient prime $\mathfrak p\mid11$,
\[
 \mathbf P(\Lambda/\mathfrak p\Lambda)\simeq\mathbf P^1(\F_{11})
\]
parametrizes lines in the relevant $\mathfrak p$-torsion local system, equivalently cyclic order-$11$ subgroups in that RM component.  The stabilizer of a line is a Borel subgroup, so the corresponding Borel level cover is a natural finite level structure on the given family.  This is a relative interpretation for these families; it is not a universal fine moduli theorem for all hypergeometric abelian varieties.

Clark--Voight construct congruence covers of a discrete triangle group
by reducing a quaternionic order; Duque-Rosero--Voight study many
low-genus Borel quotients \cite{ClarkVoight,DuqueVoight}. Here the
integral representation is first obtained from an abelian family.
Proposition~\ref{prop:CV11-endpoint} identifies its level-$11$ covers
with the Clark--Voight covers for the $\sigma_-$ constituent of
$J^-_{5,5}(3)$. The comparison respects the ordered branch values.
Corollary~\ref{cor:K-endpoint-borel} descends it to $K$. The four
covers attached to the other two families are distinguished from this
pair by their ordered trace data.

The Borel covers exist for all three families. Their construction uses
the integral homology representation and does not require a discrete
real image. The discrete constituent provides the additional
triangle-group interpretation.

\subsection{Relation with generalized Fermat geometry}

Frey abelian varieties have been studied using explicit hyperelliptic
models, Galois representations, modularity, and level lowering, with
Diophantine applications \cite{Kraus,BCDF,BCDFII,ChenKoutsianas}.
We study the varying families before specialization.
Theorem~\ref{thm:intro-kraus} relates $J^-_{5,5}(3)$ to the Frey
Jacobian in that arithmetic construction, with its quadratic twist
and change of RM embedding.

The isogeny alone does not yield level lowering, uniform big image for individual arithmetic specializations, or a generalized Fermat theorem.  Below, the geometric residual image and the cyclotomic determinant do yield a full generic arithmetic image and marked $K$-models of the torsion-line covers.  They do not determine the exceptional specialized fibres or force rational points on those covers.

\subsection{Hypergeometric context and a companion family}

The cyclic-cover description belongs to a broader theory of hypergeometric abelian varieties.  Archinard develops primitive Jacobian realizations of hypergeometric periods, while later work studies periods, endomorphisms, Galois realizations, and finite-field hypergeometric traces \cite{Archinard,DFLST,FLRST,HoffmanTu}.  We study integral Betti monodromy and its finite level covers; comparison with $\ell$-adic realizations and finite-field Frobenius traces requires the arithmetic twists and reduction hypotheses.

The companion manuscript \emph{A Hypergeometric RM Surface: Integral
Monodromy and Level-11 Covers} \cite{TTVcompanion} studies a single
genus-two Jacobian in Darmon's $(r,\infty,\infty)$ branch. Its two
real period systems have local projective orders $(5,\infty,\infty)$;
one has image $\Delta(2,5,\infty)$. The quadratic parameter
$s=4T(1-T)$ relates its projective period maps to the Hecke Gauss
equations. In the present paper the local projective orders are
$(5,5,\infty)$, and the discrete constituent gives the genuine
index-two subgroup
\[
 \Delta(5,5,\infty)<\Delta(2,5,\infty).
\]
The families are related after parameter change, base extension, and
quadratic twisting \cite{BCDF}. Their parameter lines and ordered
level covers nevertheless remain distinct. In the companion paper the
level-$11$ covers are normalizations of fibre products of Hecke covers
under $s=4T(1-T)$. Here the two covers of $J^-_{5,5}(3)$ are themselves
the corresponding triangular modular curves. Section~\ref{sec:monodromy}
explains the different quadratic maps underlying these comparisons.

\subsection{Background, proofs, and organization}

Several tools used below are standard.  Rational group-algebra idempotents describe the old and primitive factors of Jacobians with finite group actions \cite{CaroccaRodriguez}; specialized Gassner representations describe the monodromy of cyclic covers \cite{Venkataramana}, with the equal-exponent/full-braid Burau setting treated in \cite{VenkataramanaAnnals}.  Related viewpoints appear in \cite{KapovichMillson,McMullen}.  Section~\ref{sec:primitive} fixes the cyclic-cover conventions,
identifies the primitive factors needed later, and treats the index
$j=r-2$ directly.

The proof of Theorem~\ref{thm:intro-kraus} occupies Section~\ref{sec:r5kraus}, especially Subsections~\ref{sec:kraus-quotient} and \ref{sec:kraus-bcdf}.  The endpoint cyclic cover is first reduced to an explicit genus-two quotient, and an algebraic correspondence on the primitive part then gives the isogeny; Galois descent supplies the constant twist, while an explicit scaling gives the factor $1-t$ and the RM relabeling.  The marked invariant in the final subsection of Section~\ref{sec:r5kraus} separates the endpoint from the other minus family.

Section~\ref{sec:monodromy} starts from the same rank-two Gassner representations and puts them into a $K$-rational form.  Subsections~\ref{sec:nondiscrete-images} and \ref{sec:discrete-image} prove Proposition~\ref{thm:intro-real-images}: five cases contain an infinite-order elliptic element, while the remaining case is identified with the index-two group $\Delta(5,5,\infty)<\Delta(2,5,\infty)$.  Theorem~\ref{thm:intro-real-lattice} is proved in Subsection~\ref{sec:lattices}; the minus lattices follow from the local order calculation, whereas the plus case requires the integral action of the quotient involution and a local descent at the ramified prime above $5$.

Finally, Section~\ref{sec:level11} proves Theorem~\ref{thm:intro-level11}.  Subsection~\ref{sec:residual11} proves full residual image at the two primes above $11$, Subsection~\ref{sec:degree12} constructs the degree-$12$ Borel covers and identifies the endpoint pair with the corresponding Clark--Voight triangular modular curves, and Subsections~\ref{sec:borel-triangular}--\ref{sec:marked11} compare the six marked levels.  The abstract finite-group generation and Nielsen invariants used there have antecedents in \cite{McCulloughWanderley,FairbairnMagaardParker}; what matters here is their attachment to the three geometric families and to the two coefficient primes.

The hyperelliptic Frey family and its RM action are part of the existing
arithmetic construction; the expanded version of \cite{BCDF} supplies
the model over $U$ used here. We additionally determine the generic
arithmetic image at level $11$ and descend the ordered Borel covers
to $K$. Explicit Belyi equations and comparisons with a specified
abelian moduli problem, including its integral lattices and
polarizations, remain separate questions.

\section{Cyclic covers, primitive factors, and Gassner monodromy}\label{sec:primitive}

The later comparison with Frey models and torsion levels is meaningful only after the correct abelian factor has been isolated inside the cyclic-cover Jacobian.  Darmon's quotient construction supplies the relevant old factors, while rational idempotents and specialized Gassner representations provide standard tools for expressing their complements \cite{CaroccaRodriguez,Venkataramana}.  We use these tools to fix one explicit set of conventions for Darmon's systems, including the rank-one Kummer twist coming from the scalar coefficient, which distinguishes the geometric linear monodromy from its projectivization.  We also treat the endpoint $j=r-2$ directly, because the commonly quoted coprime-exponent statement does not apply there.  The section thereby fixes the precise geometric factor and marking needed for the $r=5$ hyperelliptic comparison, the real-image calculation, and the integral lattice analysis that follow.

Before introducing the group-algebra idempotents, it is useful to keep the geometry visible. The branch points move on the $x$-line while the parameter $t$ lives on a three-punctured sphere. The deck action organizes the Jacobian into cyclotomic character pieces, whereas the involution $\tau$ pairs opposite characters. The real-multiplication surface is obtained only after both operations have been taken into account. This is the reason we first work on the full cyclic cover and only later descend to the rank-two constituent.

\subsection{Darmon's cyclic covers and branch data}

Let $r$ be an odd prime and put
\[
 F=\Q(t),\qquad u=\frac{t}{t-1},\qquad
 L=\Q(\zeta_r)=\Q(\zeta_{2r}),\qquad K=L^+.
\]
For odd $j$ define
\[
 f_j(x)=u^2x^{j-2}\left(\frac{x-1}{x-u}\right)^{j+2}.
\]
We use the smooth projective normalizations of
\[
 X^-:y^{2r}=f_j(x),\qquad X^+:y^r=f_j(x),\qquad E:z^2=f_j(x).
\]
The maps relevant to Darmon's construction are
\[
 \pi:X^-\longrightarrow E,\qquad (x,y)\longmapsto(x,y^r),
\]
\[
 \pi_r:X^-\longrightarrow X^+,\qquad (x,y)\longmapsto(x,y^2),
\]
and all three curves carry the involution
\[
 \tau(x,y)=\left(\frac ux,\frac1y\right).
\]
We write
\[
 C^\pm=X^\pm/\tau,\qquad E'=E/\tau,
\]
\[
 A^\pm=\Jac(X^\pm),\qquad B^\pm=\Jac(C^\pm).
\]
Darmon's definitions are
\[
 J^+_{r,r}(j)=B^+,
\]
\[
 J^-_{r,r}(j)=B^-\Big/\left(\pi^*\Jac(E')+\pi_r^*B^+\right).
\]
The minus range is $j=1,3,\dots,r-2$, while the plus range is $j=1,3,\dots,r-4$ \cite[Prop.~1.12]{Darmon2000}.

Order the branch points as $(0,1,u,\infty)$. The Kummer exponent vectors are
\[
 \mathbf k^-(j)=(j-2,j+2,2r-j-2,2r-j+2)\pmod{2r}
\]
and
\[
 \mathbf k^+(j)=(j-2,j+2,r-j-2,r-j+2)\pmod r.
\]
These are obtained by taking the exponents at $0,1,u$ from the displayed equation and choosing the fourth exponent so that the total divisor degree is zero modulo the cyclic degree. If a local exponent is $k$ for a connected Kummer cover of degree $d$, then the ramification index is $d/\gcd(d,k)$, and Riemann--Hurwitz gives
\[
 2g(X)-2=-2d+\sum_{\nu=1}^4\bigl(d-\gcd(d,k_\nu)\bigr).
\tag{2.3}
\]
This elementary formula is enough to see the generic/endpoint dichotomy before any Jacobian decomposition is invoked. For $j\le r-4$ all four exponents are prime to the relevant degree. At the endpoint $j=r-2$ the minus vector becomes
\[
 (r-4,r,r,r+4)\pmod{2r},
\]
while the corresponding plus cover has only two ramified points and genus zero. The Riemann--Hurwitz calculation gives
\[
 g(X^-)=2r-1\quad(j\le r-4),\qquad g(X^-)=r\quad(j=r-2),
\]
\[
 g(X^+)=r-1\quad(j\le r-4).
\]
Darmon's fixed-point analysis of $\tau$ then yields the expected quotient dimensions. It is useful to record the geometry in a form that will later explain the old/new decomposition:
\[
\begin{array}{c|c|c|c|c}
\text{family}&j&\text{ramification indices}&g(X)&g(X/\tau)\\ \hline
X^-&1,3,\ldots,r-4&(2r,2r,2r,2r)&2r-1&r\\
X^-&r-2&(2r,2,2,2r)&r&(r+1)/2\\
X^+&1,3,\ldots,r-4&(r,r,r,r)&r-1&(r-1)/2.
\end{array}
\]
The omitted plus endpoint has ramification indices $(r,1,1,r)$ and genus zero. In particular, the exceptional minus family loses precisely the old order-$r$ contribution, rather than acquiring an additional factor.

We shall also use explicitly how $\tau$ exchanges cyclic characters. If $d$ is the degree of the cyclic cover and a differential in the $m$-eigenspace is written as $h(x)\,dx/y^m$, then
\[
 \tau^*\left(h(x)\frac{dx}{y^m}\right)
 =-\frac{u}{x^2}h(u/x)y^m\,dx
 =-\frac{u f_j(x)}{x^2}h(u/x)\frac{dx}{y^{d-m}}.
\tag{2.1}
\]
Thus $\tau$ interchanges the $m$ and $d-m$ eigenspaces and also their holomorphic parts. This observation is elementary, but it is important: an individual non-real character of the cyclic group is not by itself the rational hypergeometric abelian variety. The rational factor is obtained by summing a full Galois orbit and then taking the $\tau$-invariant part.

\subsection{Rational idempotents and primitive factors}

The branch calculation tells us which character orders can occur; the next task is to express Darmon's quotient using rational group-algebra idempotents. This is standard for Jacobians with a finite group action; see, for example, \cite{CaroccaRodriguez}. We record the idempotents explicitly because the later monodromy and integral calculations require a fixed identification of the relevant factors. The result is a description of Darmon's old-factor quotient in the present conventions.

Choose over $L(t)$ a deck generator
\[
 \delta(x,y)=(x,\zeta_{2r}y)
\]
on $X^-$. In the rational endomorphism algebra of $A^-$ define
\[
 P_1=\frac1{2r}\sum_{a=0}^{2r-1}\delta^a,
\quad
 P_2=\frac1r\sum_{a=0}^{r-1}\delta^{2a},
\quad
 P_r=\frac{1+\delta^r}{2},
\]
and
\[
 E_2=P_2-P_1,
\quad E_r=P_r-P_1,
\quad E_{2r}=1-P_2-P_r+P_1.
\]
Since the quotient by the full deck group is $\mathbf P^1$, the trivial idempotent $P_1$ acts as zero on the Jacobian. The three idempotents $E_2,E_r,E_{2r}$ are orthogonal and sum to the identity. Moreover they descend to $F$: Galois conjugation permutes the deck generators through $\delta\mapsto\delta^a$ with $a\in(\Z/2r\Z)^\times$, hence preserves each subgroup average.

\begin{proposition}\label{prop:old-new}
Let $A_e=\im(E_e)$ in the isogeny category. Then
\[
 A_2\sim_F\Jac(E),\qquad
 A_r\sim_F A^+,
\]
and
\[
 A^-\sim_F A_2\times A_r\times A_{2r}.
\]
Thus $A_{2r}$ is the exact-order-$2r$ primitive factor of $A^-$.
\end{proposition}

\begin{proof}
The quotient-map identities give
\[
 P_2=r^{-1}\pi^*\pi_*,\qquad
 P_r=2^{-1}\pi_r^*\pi_{r*}.
\]
Hence their nontrivial images are isogenous to the Jacobians of the corresponding quotient covers. Since the only nontrivial divisors of $2r$ are $2,r,2r$, there are no additional character-order factors.
\end{proof}

The involution satisfies
\[
 \tau\delta\tau=\delta^{-1}.
\]
Consequently it commutes with all rational subgroup averages and therefore with the $E_e$. Put
\[
 T_+=\frac{1+\tau}{2}.
\]
Pullback from $C^-$ identifies $B^-$ up to isogeny with the $T_+$-part of $A^-$. The two displayed maps $\pi$ and $\pi_r$ commute with $\tau$, so their images inside $B^-$ are exactly the $T_+E_2$ and $T_+E_r$ factors.

\begin{proposition}\label{thm:primitive-general}
For every odd prime $r$ and every odd $j$ in Darmon's minus range,
\[
 \boxed{J^-_{r,r}(j)\sim_F\im(T_+E_{2r}).}
\]
For every odd $j\le r-4$,
\[
 \boxed{J^+_{r,r}(j)\sim_F\im\left(T_+:A^+\to A^+\right),}
\]
and $A^+$ is entirely primitive of order $r$. In each sign the full primitive factor has dimension $r-1$, while the $\tau+$ hypergeometric abelian variety has dimension $(r-1)/2$.
\end{proposition}

\begin{proof}
The minus decomposition of Proposition~\ref{prop:old-new} is $\tau$-stable. Hence
\[
 B^-\sim_F\Jac(E')\times B^+\times \im(T_+E_{2r}).
\]
Darmon's quotient by the first two factors therefore induces an isogeny from $\im(T_+E_{2r})$ to $J^-_{r,r}(j)$. On the degree-$r$ cover the only nontrivial character order is $r$, so the entire Jacobian is the primitive factor; taking the $\tau+$ part gives the plus statement.
\end{proof}

\begin{remark}
Over $L(t)$ one can also compare the two $\tau$-eigenspaces. On the primitive factor the correspondence $\sigma-\sigma^{-1}$, for the relevant order-$r$ deck operator $\sigma$, anti-commutes with $\tau$ and is invertible on primitive characters. Thus the two halves are isogenous over $L(t)$ and the full primitive factor is isogenous to the square of the hypergeometric abelian variety. We shall not need this refinement below.
\end{remark}

\subsection{The exceptional endpoint}

At $j=r-2$ two branch exponents are divisible by $r$, so it would be incorrect simply to cite the generic coprime Gassner theorem. At this endpoint the geometry simplifies: the order-$r$ quotient has genus zero. We therefore inspect the deck characters directly and recover the primitive order-$2r$ piece without any appeal to the missing coprimality hypothesis.

Let $V_m$ be the eigenspace in $H^1(X^-,\C)$ for deck pullback with eigenvalue $\zeta_{2r}^{-m}$. If $s_m$ is the number of nontrivial local character multipliers, then the standard Kummer differential calculation gives
\[
 \dim H^{1,0}(X^-)_m=-1+\sum_{\nu=1}^4\left\{\frac{mk_\nu}{2r}\right\}
 =-1+\frac{s_m}{2},
\]
so that
\[
 \dim V_m=s_m-2.
\]
At $j=r-2$ the branch vector is
\[
 (r-4,r,r,r+4).
\]
If $m$ is even and nonzero, the two middle multipliers are trivial and the two outer ones are not; thus $s_m=2$ and $V_m=0$. This is the disappearance of the order-$r$ old factor. If $m$ is odd and primitive, all four multipliers are nontrivial, hence $\dim V_m=2$ and its Hodge numbers are $(1,1)$. Finally $m=r$ gives the order-$2$ elliptic factor.

It follows that
\[
 A^-\sim_F\Jac(E)\times A_{2r},
\qquad
 B^-\sim_F\Jac(E')\times A_{2r}^{\tau+},
\]
and therefore
\[
 J^-_{r,r}(r-2)\sim_F A_{2r}^{\tau+}.
\]
This proves the endpoint clause of Proposition~\ref{thm:primitive-general} without applying a coprimality-dependent theorem.

\subsection{Rank-two Gassner realization}

We now pass from the algebraic factor to its moving rank-two local system. This step explains why a rank-three reduced Gassner specialization produces a rank-two local system on compact homology and also records the scalar Kummer factor needed for the geometric linear monodromy. Keeping that scalar is essential later when we compare linear parity, even though it disappears after projectivization.

Venkataramana's Annals paper \cite{VenkataramanaAnnals} treats equal branch exponents using the Burau representation of the full braid group. His Inventiones paper \cite{Venkataramana} treats unequal exponents using the Gassner representation of the pure braid group, and is the relevant reference here. Their higher-rank arithmeticity results require a stable range such as $n\ge 2d$. Here the affine Gassner parameter is $n=3$, so those theorems do not determine discreteness or arithmeticity of our rank-two images. The real-image calculation in Section~\ref{sec:monodromy} is therefore a separate low-dimensional problem.

Fix a primitive character $m$ of the relevant cyclic degree $d$. Its local multipliers are
\[
 \mathbf t_m=
 \left(
 \zeta_d^{m(j-2)},
 \zeta_d^{m(j+2)},
 \zeta_d^{-m(j+2)},
 \zeta_d^{m(2-j)}
 \right),
\]
and their product is $1$. For the primitive factors under consideration all four entries are nontrivial, including the minus endpoint.

After moving the four branch points into an affine chart, the reduced specialized Gassner module has rank three. The condition $\prod t_i=1$ produces the peripheral invariant line. Compactification kills precisely this line, leaving a rank-two character module on compact homology. In the generic coprime cases this is the compact specialization described by Venkataramana \cite{Venkataramana}; the same compact chain argument remains valid for the primitive endpoint characters because their four local multipliers are nontrivial. For the vanished order-$r$ endpoint characters, by contrast, two multipliers are trivial and the compact $H^1$ is zero.

Here is a direct way to see the rank drop that will also be useful integrally. A rank-one local system on the affine line punctured at four finite branch points has a cellular model with four lifted meridians and one relation coming from the boundary map, hence first homology of rank three over the character field. Because the product of the four local multipliers is $1$, the loop around affine infinity defines a nonzero class. Its coefficients are successive partial products of the multipliers, so it is nonzero whenever the four individual multipliers are nontrivial. Filling the point at infinity kills this class. Filling the branch points creates no further relation in the primitive character because none of the local multipliers equals $1$. Thus compact homology has rank exactly two. This argument does not require the branch exponents themselves to be units modulo the full cyclic degree; it uses only the nontriviality of the specialized primitive multipliers.

The reduced Gassner action preserves this peripheral line, so the local system on the homology of the compactified cover is obtained by a quotient representation, not merely by deleting an eigenvalue after the fact. The specialized skew-Hermitian form has the same peripheral radical and descends to a nondegenerate form on the rank-two quotient. This is the character local system on compact homology that enters the weight-one Hodge structure of the hypergeometric abelian variety.

The involution exchanges opposite characters. If $V_m$ denotes one character local system on compact homology, then
\[
 (V_m\oplus V_{-m})^{\tau+}
 =\{v+\tau^*v:v\in V_m\}\simeq V_m.
\]
Thus each real-multiplication constituent is identified with the rank-two compact-homology quotient of a specialized Gassner representation.  This is not merely a statement about the open cyclic cover: the rank-three reduced Gassner specialization has a one-dimensional peripheral radical, and quotienting by that line is exactly what fills the auxiliary point at infinity.

There is one further normalization that becomes essential when traces rather than projective classes are compared. Move the four branch points to finite points by
\[
 x=2+\frac1z.
\]
On the open set $u\ne2$ the equation can be written in the exact form
\[
 y^d=C_j(t)(z+1/2)^{j-2}(z+1)^{j+2}
 \left(z-\frac1{u-2}\right)^{-j-2}z^{2-j},
\]
where
\[
 C_j(t)=\frac{u^2 2^{j-2}}{(2-u)^{j+2}}.
\tag{2.2}
\]
If $\eta^d=C_j(t)$ and $\widetilde y=y/\eta$, then on the $m$-character
\[
 \frac{dx}{y^m}=\eta^{-m}\frac{dx}{\widetilde y^m}.
\]
Thus the usual unit-leading-coefficient Gassner family differs from the geometric linear local system by a rank-one Kummer character. The point $u=2$ is only an apparent singularity introduced by this coordinate choice and disappears when the configuration braid and scalar factor are recombined. This is why the projective monodromy is easier to see than the precise linear lift, and it explains the central-twist discussion in Section~\ref{sec:monodromy}.

The rational abelian factors are thus described by idempotents selecting
characters of the specified orders. Each RM constituent is the rank-two
quotient of the specialized Gassner representation on the homology of
the compactified cover. The direct argument at $j=r-2$ supplies this
description without a coprimality hypothesis. We next compare that
family with the Kraus--BCDF Jacobian, and then identify the integral
homology lattices inside the rational representations.

\section{The three r=5 families and the Kraus correspondence}\label{sec:r5kraus}

For $r=5$, Darmon's $(r,r,p)$ construction gives several abelian-surface
families with local projective orders $(5,5,\infty)$. Dimension, RM
field, and parity therefore do not by themselves identify a particular
Frey model. Kraus constructs the relevant genus-two RM family, and
Billerey--Chen--Dieulefait--Freitas give a model over $U$ and place it
among Darmon's odd representations \cite{Kraus,BCDF}. We identify its
index $j$ by constructing a genus-two quotient at $j=3$ and an isogeny
over $K(t)$, with its quadratic twist and nontrivial change of RM
embedding. This comparison attaches the later monodromy and level
calculations to the Frey Jacobian.

We first record the cyclic degrees, branch vectors and linear parity of the three families. The endpoint requires separate treatment because its branch exponents are not all coprime to the cyclic degree.

\subsection{The three r=5 families}

From now on
\[
 K=\Q(\sqrt5),\qquad R=\mathcal O_K=\Z[\alpha],\qquad
 \alpha=\frac{\sqrt5-1}{2}.
\]
There are exactly three nonzero Darmon hypergeometric abelian varieties. Their cyclic-cover data are summarized in Table~\ref{tab:r5inventory}.

\begin{table}[ht]
\centering
\begin{tabular}{c c c c c}
\toprule
hypergeometric abelian variety & degree $d$ & branch vector & primitive $m$ & dimension \\
\midrule
$J^-_{5,5}(1)$ & $10$ & $(9,3,7,1)$ & $1,3$ & $2$\\
$J^-_{5,5}(3)$ & $10$ & $(1,5,5,9)$ & $1,3$ & $2$\\
$J^+_{5,5}(1)$ & $5$ & $(4,3,2,1)$ & $1,2$ & $2$\\
\bottomrule
\end{tabular}
\caption{The three $r=5$ Darmon families. Branch points are ordered $(0,1,u,\infty)$.}
\label{tab:r5inventory}
\end{table}

The first minus row is generic. The second is the endpoint $j=3$ and has no order-$5$ old factor. The plus row is generic. Each primitive inverse pair determines one of the two real embeddings of $K$.

For comparison with the Frey-variety literature, recall the hyperelliptic curve used by Kraus. If
\[
 h_r(T)=\prod_{m=1}^{(r-1)/2}
 \left(T-(\zeta_r^m+\zeta_r^{-m})\right),
\]
then a solution of $a^r+b^r=c^p$ leads to
\[
 C_r(a,b):\qquad
 y^2=(ab)^{(r-1)/2}x\,h_r\!\left(\frac{x^2}{ab}+2\right)+b^r-a^r.
\]
Kraus proves that this curve has genus $(r-1)/2$ and real multiplication by $\Q(\zeta_r)^+$. The later three-point model of \cite{BCDF} is particularly well suited to comparison with Darmon's parameter $t$. Their construction proves that the resulting Jacobian gives an odd Frey system and hence belongs on Darmon's minus side. In the sources and conventions compared here, the complete assignment to a Darmon index together with the twist and RM generator is not stated in the form needed for our later normalization. The next calculation supplies that explicit comparison.

\subsection{A genus-two quotient of the endpoint cover}\label{sec:kraus-quotient}

The endpoint degree-ten cover has genus five before taking the relevant quotient.  We now expose a second involution whose quotient has genus two and keep track of the fields of definition throughout.  In this subsection all affine equations are understood through the smooth projective normalizations of their function fields.  Put
\[
 F=\Q(t),\qquad K=\Q(\sqrt5),\qquad L=\Q(\zeta_5),\qquad u=\frac{t}{t-1}.
\]
The endpoint minus model is
\[
 y^{10}=u^2x\left(\frac{x-1}{x-u}\right)^5.
\]
Define over $F$
\[
 v=y^2\frac{x-u}{x-1},\qquad
 w=y(v^5-u^3).
\]
Then $v^5=u^2x$, and the function field is birational over $F$ to
\[
 H:\qquad w^2=v(v^5-u^2)(v^5-u^3),
\]
with inverse functions
\[
 x=\frac{v^5}{u^2},\qquad y=\frac{w}{v^5-u^3}.
\]
The polynomial on the right has degree $11$ and is generically squarefree, so the smooth projective model of $H$ has genus five.

In these coordinates
\[
 \tau(v,w)=\left(\frac uv,-\frac{u^3w}{v^6}\right),
\]
while the hyperelliptic involution is $h(v,w)=(v,-w)$.  Put $\iota=h\tau$.  The functions
\[
 \xi=v+\frac uv,
 \qquad
 z=\frac{w}{v^3}
\]
are $\iota$-invariant and satisfy
\[
 Q_u:\qquad
 z^2=F_u(\xi),\qquad
 F_u(Z)=Z^5-5uZ^3+5u^2Z-u^2(1+u).
\]
Since $v$ satisfies
\[
 v^2-\xi v+u=0,
\]
the invariant subfield has degree two and is exactly the fixed field of $\iota$.  Thus $Q_u=H/\iota$ as a smooth projective quotient, not merely as a rational image; in particular $Q_u$ has genus two.  The quotient maps
\[
 q_+:H\longrightarrow H/\tau=C^-,\qquad
 q_-:H\longrightarrow H/\iota=Q_u
\]
are defined over $F$.

Let $\zeta=\zeta_5$ and let $\sigma(x,y)=(x,\zeta y)$ on the original cover.  On $H$ one has
\[
 \sigma(v,w)=(\zeta^2v,\zeta w).
\]
Thus $\sigma$, and hence the correspondence below, is defined over $L(t)$.
On Jacobians the symbol $\sigma$ denotes the covariant action
$\sigma_*$; its adjoint for the canonical polarization is
$\sigma^{-1}$.  Put
\[
 D=\sigma-\sigma^{-1}.
\]
Let $B^-=\Jac(C^-)$ and let $O\subset B^-$ be the connected old elliptic subvariety coming from the Legendre quotient.  The commuting quotient diagram shows that $q_+^*O$ lies in the pullback of the order-two quotient of $H$.  The order-five deck transformation $\sigma$ acts trivially on that quotient, so
\[
 Dq_+^*|_O=0.
\]
Consequently
\[
 \widetilde\Phi=q_{-*}Dq_+^*:B^-\longrightarrow\Jac(Q_u)
\]
vanishes on $O$ and factors, by the universal property of the abelian quotient, through
\[
 B^-/O=J^-_{5,5}(3).
\]
There is no additional order-five old Jacobian at this endpoint: the degree-five quotient has genus zero.

The resulting homomorphism is an isogeny for a rank reason, not from simplicity.  On every primitive eigenspace of the order-ten part, $D$ has eigenvalue $\zeta^a-\zeta^{-a}\ne0$.  Moreover $D\tau=-\tau D$, so $D$ carries the full primitive $\tau$-invariant subspace isomorphically onto the primitive $\tau$-anti-invariant subspace.  Since $h$ acts as $-1$ on $H_1(H,\Q)$, the latter is the $\iota$-invariant primitive subspace.  Pullback by $q_+$ identifies the rational homology of $J^-_{5,5}(3)$ with the first space, while on the second
\[
 q_-^*q_{-*}=1+\iota=2.
\]
Hence the induced map on rational first homology has rank four.  Both source and target have dimension two, so the factored homomorphism
\[
 \Phi:J^-_{5,5}(3)\longrightarrow\Jac(Q_u)
\]
is an isogeny over $L(t)$.

We next descend this isogeny to $K(t)$.  For a hyperelliptic curve
\[
 Q:z^2=F(\xi)
\]
and $a\in K(t)^\times$, write
\[
 Q^{[a]}:Z^2=aF(\xi)
\]
for our quadratic-twist convention.  Put
\[
 \theta=\zeta+\zeta^{-1}=\alpha,
 \qquad
 \eta=\zeta-\zeta^{-1},
 \qquad
 b=\eta^2=\alpha^2-4=-3-\alpha=-\frac{5+\sqrt5}{2}.
\]
Let $c$ be the nontrivial element of $\Gal(L/K)$.  It fixes $H$, both quotient maps, and the involutions, but sends $\sigma$ to $\sigma^{-1}$; hence
\[
 c(\Phi)=-\Phi.
\]
Over $L(t)$ the map
\[
 j_b:Q_u^{[b]}\longrightarrow Q_u,
 \qquad
 (\xi,Z)\longmapsto(\xi,Z/\eta)
\]
is an isomorphism.  Since $c(\eta)=-\eta$, one has $c(j_b)=h\circ j_b$, where $h$ is the hyperelliptic involution.  On Jacobians this says $c(j_{b*})=-j_{b*}$.  Therefore
\[
 \Phi_b=j_{b*}^{-1}\Phi:
 J^-_{5,5}(3)\longrightarrow\Jac(Q_u^{[b]})
\]
is $c$-invariant:
\[
 c(\Phi_b)=(-j_{b*})^{-1}(-\Phi)=\Phi_b.
\]
Galois descent of morphisms gives the $K(t)$-isogeny
\[
 J^-_{5,5}(3)\sim_{K(t)}\Jac(Q_u^{[b]}).
\]

\subsection{Comparison with the BCDF three-point model}\label{sec:kraus-bcdf}

The quotient just obtained is already a genus-two hyperelliptic family.  We now put it into the same parameter and normalization as the expanded BCDF construction.  Two twists must remain distinct: the constant factor $b$, which descends the correspondence from $L$ to $K$, and the variable factor $1-t$, which changes the three-point model.

The genus-two family of \cite{BCDF} is
\[
 C'_5(t):\qquad
 Y^2=F_{C'}(X):=X^5+5qX^3+5q^2X+q^2(2t-1),
 \qquad q=t(1-t).
\]

\begin{remark}[The three-point parameter convention]\label{rem:ck-table-sign}
The polynomial above follows the expanded BCDF convention
\cite[Lemma~2.36]{BCDF}.  The version printed in
\cite[Table~1.4]{ChenKoutsianas} has constant term
$-q^2(2t-1)$ rather than $+q^2(2t-1)$.  With the same
$q=t(1-t)$ it is $F_{C'}(X,1-t)$; if the parameter $t$
is kept fixed, the two equations are related instead by
$F_{\mathrm{CK}}(X,t)=-F_{C'}(-X,t)$, requiring the
corresponding quadratic twist by $-1$.  We retain the BCDF
convention throughout the statement of Theorem~\ref{thm:kraus-main};
no sign is absorbed into the constant twist $b$.
\end{remark}

Set $d=1-t$.  Since $u=t/(t-1)$, direct substitution gives the identity
\[
 F_{C'}(d\xi)=d^5F_u(\xi).
\]
It follows that
\[
 Q_u^{[d]}\longrightarrow C'_5(t),
 \qquad
 (\xi,Z)\longmapsto(d\xi,d^2Z),
\]
is an $F$-isomorphism.  More directly, the constant-twisted model appearing in the descent is carried over $K(t)$ to the precise twist used below by
\[
 Q_u^{[b]}\longrightarrow (C'_5(t))^{[bd]},
 \qquad
 (\xi,Z)\longmapsto(d\xi,d^3Z),
\]
because
\[
 (d^3Z)^2=bd^6F_u(\xi)=bd\,F_{C'}(d\xi).
\]
Thus the coordinate change gives the stated twist $b(1-t)$.

\begin{theorem}\label{thm:kraus-main}
There is an isogeny
\[
 \boxed{
 J^-_{5,5}(3)\sim_{K(t)}\Jac(C'_5(t))^{[b(1-t)]}.}
\]
Let
\[
 i_D:K\hookrightarrow\operatorname{End}^0\!\left(J^-_{5,5}(3)\right),
 \qquad
 i_B:K\hookrightarrow\operatorname{End}^0\!\left(\Jac(C'_5(t))\right)
\]
denote the Darmon and BCDF real-multiplication actions, respectively.  If $a\mapsto\bar a$ denotes the nontrivial automorphism of the abstract RM field $K$, the isogeny $\Psi$ above may be chosen so that
\[
 \Psi\,i_D(a)=i_B(\bar a)\,\Psi\qquad(a\in K).
\]
In particular the two named real embeddings are exchanged.
\end{theorem}

\begin{proof}
The algebraic isogeny and its field of definition were proved above; it remains to compare the two RM labels on that same map.  Let
\[
 S=\sigma+\sigma^{-1}
\]
on the primitive part.  The operator $D=\sigma-\sigma^{-1}$ commutes with $S$, and $S$ commutes with the quotient involutions, so the correspondence defining $\Phi$ intertwines the descended $S$-actions.

For the BCDF auxiliary model, Proposition~2.25 and equations (2.27)--(2.29) of the expanded version of \cite{BCDF} use the order-five rotation
\[
 \beta:(V,W)\longmapsto(\zeta V,\zeta^3W).
\]
With $d=1-t$, the scalar coordinate change
\[
 V=dv,\qquad W=d^5\sqrt d\,w
\]
over $F(\sqrt d)$ commutes with the rotations and carries $\beta$ to $\sigma^3$.  On the quotient it induces precisely the scaling $X=d\xi$ used above; the transported RM convention on the three-point family is recorded in \cite[(2.39)]{BCDF}.  Therefore, on the primitive factor,
\[
 \beta+\beta^{-1}
 =\sigma^3+\sigma^{-3}
 =S^2-2.
\]
If $\theta=\alpha$ is the abstract generator of $K$, this gives
\[
 \theta\longmapsto\theta^2-2=-1-\theta,
\]
which is the nontrivial automorphism of $K$.  The constant and variable quadratic twists are central and hence preserve this intertwining relation.  This proves the displayed identity for $i_D$ and $i_B$.
\end{proof}

\begin{corollary}[Degree and integral comparison]
\label{cor:reverse20}
For the isogeny $\Psi$ constructed in Theorem~\ref{thm:kraus-main},
there is a homomorphism
$\Xi:\Jac(C'_5(t))^{[b(1-t)]}\to J^-_{5,5}(3)$ over $K(t)$ such that
\[
 \Xi\Psi=[20]_{J^-_{5,5}(3)},\qquad
 \Psi\Xi=[20]_{\Jac(C'_5(t))^{[b(1-t)]}}.
\]
More precisely, $\Psi=[2]V'$ for an isomorphism
\[
 V':J^-_{5,5}(3)\xrightarrow{\sim}\Jac(C'_5(t))^{[b(1-t)]}
\]
over $K(t)$, with the same RM relabeling as $\Psi$.
Consequently $\deg\Psi=16$ and $\ker\Psi=J^-_{5,5}(3)[2]$.
For every coefficient prime $\mathfrak l$ of residue
characteristic different from $2$, the map $\Psi$ induces an
integral Tate-module isomorphism and an isomorphism of
residual Galois representations.  The isomorphism $V'$ gives
these comparisons at every coefficient prime.  With $\kappa$ the nontrivial
$K/\Q$-automorphism, these comparisons exchange $\mathfrak l$
with $\kappa(\mathfrak l)$ on the BCDF side and retain the twist
character $\chi_{b(1-t)}$.
\end{corollary}
\begin{proof}
Use the degree-two maps $q_+,q_-$ and quotient
$\pi:\Jac(C^-)\to J=J^-_{5,5}(3)$ from Section~\ref{sec:kraus-quotient}.
Write $D=\sigma-\sigma^{-1}$ and $S=\sigma+\sigma^{-1}$.
The forward morphism over $L(t)$ satisfies
$\Phi\pi=q_{-*}Dq_+^*$.  Its reverse candidate
$\Theta=\pi q_{+*}Dq_-^*$ is an actual homomorphism;
no rational projector is involved.  The integral symmetric
correspondence $T=q_{+*}\sigma q_+^*$ descends through the
old elliptic factor and induces $i_D(\theta)$ on $J$, with
$\theta=\zeta_5+\zeta_5^{-1}=\alpha$.
The degree-two norm identities, the relation
$D^2=S^2-4$, and the invariance of $Dq_+^*$ under
$\iota=h\tau$ give
\[
 \Theta\Phi\pi
 =\pi q_{+*}D(1+\iota)Dq_+^*
 =4i_D(\theta^2-4)\pi.
\]
Since $\pi$ is surjective, $\Theta\Phi=4i_D(b)$,
where $b=\theta^2-4=-3-\theta$.  Both $\Phi$ and $\Theta$
change sign under $\Gal(L/K)$; composing with the same
constant-twist isomorphism $j_b$ used above makes both
maps $K(t)$-rational.  The variable coordinate change
$Q_u^{[b]}\simeq(C'_5(t))^{[b(1-t)]}$ transports them
without changing their composition.  Finally
$\kappa(b)=-2+\theta$, $b\kappa(b)=5$, and
$i_D(\kappa(b))=T_J-[2]$ is integral.  Multiplying the
transported reverse map by $i_D(\kappa(b))$ gives
$\Xi\Psi=[20]$.  Surjectivity of $\Psi$ gives the other
composition.  In particular this already gives the integral inverse
$20^{-1}\Xi$ away from $2$ and $5$.

To determine the degree, one must use the quotient polarization,
not just these composition identities.  Put $A=\Jac(Q_u)$ and
$N=q_{-*}\sigma q_+^*:B^-\to A$.
Since $\tau\sigma\tau=\sigma^{-1}$ and $\tau$ acts as $-1$
on $q_-^*A$, one has $q_{-*}\sigma^{-1}q_+^*=-N$ and
$q_{-*}q_+^*=0$.  As $\sigma$ fixes $q_+^*O$, the integral
map $N$ kills $O$.  Thus $N=V\pi$ for a homomorphism
$V:J\to A$, and $\Phi=2V$.

The old elliptic factor has restricted polarization $[5]$.
Indeed $H/\langle\sigma\rangle$ has equation
$r^2=s(s-u^2)(s-u^3)$, with $s=v^5$, $r=wv^2$.
The induced involution
$(s,r)\mapsto(u^5/s,-u^5r/s^2)$ is a fixed-point-free
translation.  The commuting quotient diagram therefore gives
a degree-five map $e:C^-\to E'$, where $E'$ is elliptic.
Its pullback is injective: a nonzero point in its kernel,
which is killed by $5$, would give an intermediate connected
\'etale degree-five cover of $E'$, leaving a degree-one map
from the genus-three curve $C^-$ to an elliptic curve.
Hence $O=e^*\Jac(E')$ is embedded and has restricted
polarization $[5]$.
Let $\lambda_{B^-}$ be the principal polarization and put
\[
 a=\pi\lambda_{B^-}^{-1}\pi^\vee:J^\vee\longrightarrow J.
\]
The embedding $\lambda_{B^-}^{-1}\pi^\vee$ identifies
$J^\vee$ with the connected orthogonal complement of $O$.
Its intersection with $O$ is the kernel of the restricted
polarization on $O$, of order $25$.  Thus $a$ has type
$(1,5)$ and degree $25$.

All adjoints in the next calculation use canonical Jacobian
polarizations.  With $S_A=q_{-*}\sigma q_-^*$, the two norm
identities and $\iota=-\tau$ give
\[
 NN^\dagger=2-q_{-*}\sigma^2q_-^*,\qquad
 S_A^2=2+q_{-*}\sigma^2q_-^*,\qquad
 S_A^2+S_A-1=0.
\]
The last equality is the primitive order-five RM relation;
the invariant elliptic part is $\tau$-invariant and is absent
from $q_-^*A$.  Moreover $V i_D(\theta)=S_A V$.
It follows that
\[
 V a V^\vee\lambda_A=NN^\dagger=3+S_A,
 \qquad
 25(\deg V)^2=\deg(3+S_A)=N_{K/\Q}(3+\theta)^2=25.
\]
Therefore $V$ is an isomorphism.  Transport through $j_b$
and the displayed variable-twist isomorphism gives $V'$.
It descends to $K(t)$ because $2V'=\Psi$ descends and
homomorphism groups are torsion-free.  This proves the degree
and kernel assertions.  The inverse on Tate modules away
from $2$ is $\tfrac12(V')^{-1}$; the RM exchange is unchanged.
The maps and their inverses extend over the common smooth
base, so the comparisons specialize there.
\end{proof}

\begin{remark}
The isomorphism $V'$ in Corollary~\ref{cor:reverse20} gives a
principal polarization $\lambda_0=(V')^\vee\lambda_{A'}V'$ on
the quotient $J$, where $A'=\Jac(C'_5(t))^{[b(1-t)]}$.
Relative to this specified normalization,
$\Psi^\vee\lambda_{A'}\Psi=4\lambda_0$; it is not the
principal form itself.  The natural quotient polarization
$\lambda_q=5a^{-1}$ instead has type $(1,5)$ and satisfies
$\lambda_q=\lambda_0\circ i_D(2-\theta)$.
Thus no pre-existing canonical principal polarization on the
Darmon quotient has been assumed.  The independent lattice
proof below retains its original normalization.
\end{remark}

\subsection{A marked invariant distinguishing the two minus families}

A projective monodromy invariant also distinguishes the endpoint $j=3$ from the other minus family. Normalized squared traces are insensitive to scalar twists and therefore compare the two minus systems at the projective level while still remembering which finite branch point is which.

For a rank-two local monodromy matrix $T_s$ define
\[
 \kappa_s=\frac{\Tr(T_s)^2}{\det(T_s)}.
\]
This is unchanged by scalar twists. For a full twist around two colliding cyclic branch points, the ratio of the two projective eigenvalues is the product of their local character multipliers, up to inversion. Therefore
\[
 \kappa=2+\lambda+\lambda^{-1}.
\]
At $t=0$ the colliding pair is $(0,u)$, while at $t=1$ it is $(u,\infty)$. For the degree-ten minus systems and primitive representatives $m=1,3$ this gives
\[
 \kappa_0(m)=4\cos^2\frac{2\pi m}{5},
 \qquad
 \kappa_1(m)=4\cos^2\frac{\pi jm}{5}.
\]
Writing
\[
 a_-=\frac{3-\sqrt5}{2},\qquad
 a_+=\frac{3+\sqrt5}{2},
\]
we obtain
\[
\begin{array}{c|cc}
 &m=1&m=3\\ \hline
J^-_{5,5}(1)&(a_-,a_+)&(a_+,a_-)\\
J^-_{5,5}(3)&(a_-,a_-)&(a_+,a_+).
\end{array}
\]
The equality or inequality of the two marked coordinates survives inversion of character, interchange of the two real embeddings, scalar quadratic twisting, and interchange of the two finite elliptic punctures. Thus, under these marked operations, the invariant excludes the $j=1$ orbit and is consistent with the algebraic identification in Theorem~\ref{thm:kraus-main}.  It is not used to prove the descent or the isogeny itself.

\section{Exact monodromy, real projective images, and integral structures}\label{sec:monodromy}

A rational hypergeometric representation does not yet determine either its global real geometry or its finite levels.  The first question is archimedean: among the six real RM constituents, which projective representations are actually discrete triangle groups rather than merely having the same local orders?  The second is integral: which monodromy-stable $R$-lattice is the saturated lattice inherited from Jacobian homology?  The trace and Hecke-group arguments used for the real images are classical, and the cyclic-cover/Gassner model is standard; our contribution is to compare all six constituents in common $K$-coordinates and then identify the geometric lattice inside it.  The plus family requires an additional integral argument: rational monodromy leaves two possible local lattice classes at the ramified prime above $5$, and the semilinear involution descent selects one.  After this lattice is identified, reduction at a coefficient prime computes
the torsion local system of the abelian family and defines the covers
of Section~\ref{sec:level11}.

The discussion naturally has two stages. We first work over $K$ and over its two real embeddings, where the key issue is global discreteness. We then return to integral homology. Keeping these stages separate avoids a common ambiguity: a matrix representation can be perfectly well defined over $K$ while several non-homothetic $R$-lattices remain stable under it.

\subsection{Pure braids and rank-two Gassner matrices}

We fix one base point and one braid convention once and for all. The moving branch point $u=t/(t-1)$ collides successively with $0$, $\infty$, and $1$ as $t$ approaches $0$, $1$, and $\infty$. This translates the three peripheral loops in the $t$-line into explicit pure braids, so the compact Gassner formulas can be evaluated without analytic continuation or numerical approximation.

We now compute the geometric monodromy of the three $r=5$ families. Choose the base point $t_*=1/2$, so $u_*=-1$, and order the three finite branch points as $(u,0,1)$. With Venkataramana's left Artin action, let
\[
 A_{12}=s_1^2,
 \qquad
 A_{13}=s_2^{-1}s_1^2s_2.
\]
We choose positive based loops so that
\[
 \gamma_0=A_{12},\qquad
 \gamma_\infty=A_{13},\qquad
 \gamma_0\gamma_1\gamma_\infty=1.
\]
At $t=0,1,\infty$, respectively, the moving point $u$ collides with $0,\infty,1$.

For one primitive character put
\[
 q=\zeta_d,
 \qquad
 z=q^{mk_u},\qquad
 v=q^{mk_0},\qquad
 w=q^{mk_1}=z^{-1}.
\]
In the compact rank-two Gassner basis the two pure braids act as
\[
 G_0=\begin{pmatrix}zv&1-z\\0&1\end{pmatrix},
\qquad
 G_\infty=
\begin{pmatrix}
 w&1-w\\
 w-1&1-w(1-z)
\end{pmatrix}.
\]
The second matrix is nontrivially unipotent. These formulas can be obtained directly from the Artin action on the free group and Fox differentiation. In particular $A_{13}$ is computed from the chosen braid: its off-diagonal entries retain the marking of the three punctures. This marking is what later distinguishes systems having the same unordered local orders.

Retaining the coefficient $u^2$ in Darmon's Kummer equation contributes the scalar $q^{2m}$ around $u=0$. Thus the geometric homological monodromy is
\[
 T_0=q^{2m}G_0,
 \qquad
 T_\infty=G_\infty,
 \qquad
 T_1=T_0^{-1}T_\infty^{-1},
\]
and
\[
 T_0T_1T_\infty=I.
\]
For the endpoint $j=3$, the compact rank-two module is justified by the direct argument of Section~\ref{sec:primitive}, not by a coprimality hypothesis.

The determinant calculation is also instructive. Before the scalar correction, $G_0$ has determinant $zv$, whereas $G_\infty$ has determinant $zw=1$. The factor $q^{2m}$ corrects the first determinant so that the final rank-two matrices lie in $\SL_2$ for all six constituents. Thus the scalar Kummer term is not cosmetic: without it one would have the correct projective braid representation but the wrong linear determinant character.

\subsection{The K-rational form}

The monodromy matrices on the cyclic character eigenspaces are defined
over $\Q(\zeta_5)$ or $\Q(\zeta_{10})$, and the involution pairs opposite
characters. Passing to a basis adapted to that pairing descends each rank-two constituent to the real field $K=\Q(\sqrt5)$. The resulting normal form is especially convenient because the three families differ only in two parameters $x$ and $h$.

The vector $e=(1,1)^t$ is fixed by $T_\infty$. Put
\[
 P_0=[e,T_0e].
\]
A direct conjugation gives
\[
 P_0^{-1}T_0P_0=
 A(x)=\begin{pmatrix}0&-1\\1&x\end{pmatrix},
\]
\[
 P_0^{-1}T_\infty P_0=
 P(h)=\begin{pmatrix}1&h\\0&1\end{pmatrix},
\]
\[
 P_0^{-1}T_1P_0=
 B(x,h)=\begin{pmatrix}x&1-xh\\-1&h\end{pmatrix}.
\]
All determinants are $1$ and $A(x)B(x,h)P(h)=I$. For the first constituent of each family the parameters are
\[
\begin{array}{c|ccc}
 &x&h&\Tr B=x+h\\ \hline
J^-_{5,5}(1)&\alpha&1&\phi\\
J^-_{5,5}(3)&\alpha&-2\alpha&-\alpha\\
J^+_{5,5}(1)&-\phi&\sqrt5&\alpha.
\end{array}
\]
The second constituent is obtained by applying the nontrivial automorphism of $K$. The same cyclotomic change of basis $P_0$ puts all three matrices into a form defined over $K$. This point becomes important in the integral descent, where different rational conjugators can represent different lattice classes at ramified primes.

The normalized trace pairs recover the invariants of the preceding section. In particular the two minus systems remain distinct after projectivization. For reference, the six linear local systems have the following eigenvalue patterns before passing to $K$-coordinates. If $q=\zeta_d$, then
\[
\begin{array}{c|c|c}
\text{constituent}&(T_0,T_1)&T_\infty\\ \hline
J^-_{5,5}(1),m=1&(q^{\pm2},q^{\pm1})&(1,1)\\
J^-_{5,5}(1),m=3&(q^{\pm6},q^{\pm3})&(1,1)\\
J^-_{5,5}(3),m=1&(q^{\pm2},q^{\pm3})&(1,1)\\
J^-_{5,5}(3),m=3&(q^{\pm6},q^{\pm9})&(1,1)\\
J^+_{5,5}(1),m=1&(q^{\pm2},q^{\pm1})&(1,1)\\
J^+_{5,5}(1),m=2&(q^{\pm4},q^{\pm2})&(1,1).
\end{array}
\]
The infinity matrix is nontrivially unipotent in every row. Consequently all six projective representations have local orders
$(5,5,\infty)$, although their linear lifts and global images differ.

\subsection{Linear monodromy and the odd/even distinction}

Projective monodromy forgets central signs, whereas Darmon's odd/even distinction is a statement about a preferred linear lift. We therefore pause before studying real images to identify which matrices represent the geometric continuation on Jacobian homology and how the quadratic character of $1-t$ changes that lift. This removes a potential ambiguity in the order-five versus order-ten finite generators.

The matrices above describe the geometric Jacobian local systems. In the two minus cases
\[
 A^5=I,\qquad B^5=-I,
\]
while $P$ is unipotent with trace $2$. Thus the product-one lift necessarily contains an order-$10$ peripheral matrix. This is compatible with Darmon's definition of oddness: the canonical order-five lifts of the two finite projective monodromies, together with the trace-$2$ lift at infinity, have product $-I$ \cite[Def.~1.4]{Darmon2000}. Tensoring by the quadratic character of $1-t$, with local signs
\[
 (+1,-1,-1)
\]
at $(0,1,\infty)$, replaces $(A,B,P)$ by $(A,-B,-P)$. The resulting quadratic twist is the lift used for the odd Frey
representation: its matrix at infinity has trace $-2$. The projective
representation is unchanged. In the plus case $A^5=B^5=I$ and the canonical product is $I$.

This distinction is relevant to linear local systems but irrelevant to the real projective classification below.

\subsection{Rotation classes and nondiscrete images}\label{sec:nondiscrete-images}

Once the exact real matrices are known, discreteness can be tested globally. The local orders alone are insufficient: all six constituents have
projective local orders $(5,5,\infty)$. Our strategy is to find, in each nondiscrete case, a short word which is elliptic at the chosen real embedding but whose algebraic conjugate has trace outside $[-2,2]$. Such an element must be an infinite-order elliptic and therefore cannot lie in a discrete Fuchsian group.

Let $\sigma_+(\sqrt5)>0$ and $\sigma_-(\sqrt5)<0$. The two real constituents of each family preserve the split alternating form
\[
 \Omega=\begin{pmatrix}0&1\\-1&0\end{pmatrix}
\]
and therefore act orientation-preservingly on the upper half-plane. The rotation numbers of the two finite elliptic generators, measured in units of $2\pi/5$, are
\[
\begin{array}{c|cc}
 &\sigma_+&\sigma_-\\ \hline
J^-_{5,5}(1)&(3,1)&(1,3)\\
J^-_{5,5}(3)&(3,3)&(1,1)\\
J^+_{5,5}(1)&(1,2)&(3,4).
\end{array}
\]

Five of these representations are nondiscrete. We use a simple exact criterion. If an elliptic element of $\SL_2(\mathbf R)$ has finite projective order, its trace is a sum of roots of unity; every algebraic conjugate of that trace therefore lies in $[-2,2]$. Hence an elliptic element whose trace has one real conjugate outside $[-2,2]$ has infinite order. A discrete Fuchsian group cannot contain such an infinite-order elliptic element.

The following table gives such an element in each nondiscrete case.
\[
\begin{array}{c|c|c|c}
\text{representation}&\text{word}&\text{trace}&\text{conjugate trace}\\ \hline
J^-_{5,5}(1),\sigma_+&AP^{-1}&\alpha-1&-\phi-1\\
J^-_{5,5}(1),\sigma_-&AP^2&2-\phi&2+\alpha\\
J^-_{5,5}(3),\sigma_+&AP^{-1}&3\alpha&-3\phi\\
J^+_{5,5}(1),\sigma_+&A^2PAP&2-\sqrt5&2+\sqrt5\\
J^+_{5,5}(1),\sigma_-&A^2P^{-1}&2-\sqrt5&2+\sqrt5.
\end{array}
\]
In every row the indicated trace has absolute value less than $2$, while its conjugate has absolute value greater than $2$.

\subsection{The unique discrete constituent}\label{sec:discrete-image}

The remaining real constituent behaves very differently. Instead of using the trace criterion applied above, we identify its generators inside the classical Hecke group and show that they generate the expected index-two free product $C_5*C_5$. This gives faithfulness, not merely discreteness.

It remains to treat $J^-_{5,5}(3)$ under $\sigma_-$. Here
\[
 x=-\phi,
 \qquad
 h=2\phi.
\]
Set
\[
 S=\begin{pmatrix}0&-1\\1&0\end{pmatrix},
 \qquad
 T=\begin{pmatrix}1&\phi\\0&1\end{pmatrix}.
\]
Since $\phi=2\cos(\pi/5)$, the projective group generated by $S$ and $T$ is the classical Hecke group of parameter $5$. Our generators satisfy
\[
 A=ST^{-1},\qquad P=T^2.
\]
Projectively, the Hecke group is $C_2*C_5$. The subgroup generated by $A$ and $P$, equivalently by $A$ and $SAS$, is the kernel of the homomorphism
\[
 C_2*C_5\longrightarrow C_2
\]
that is nontrivial on the $C_2$ factor and trivial on the $C_5$ factor. Hence it is
\[
 C_5*C_5\cong\Delta(5,5,\infty).
\]
Moreover
\[
 B=A^{-1}(SAS)A,
\]
so the marked finite elliptic generators give a Nielsen-equivalent free-product generating pair. This proves Proposition~\ref{thm:intro-real-images}. We use here only the standard Hecke-group presentation; see, for example, \cite{SahinBizimCangul}.

\begin{proposition}[Endpoint restriction of the Hecke periods]
\label{prop:endpoint-hecke-periods}
For the discrete real constituent of $J^-_{5,5}(3)$, put
\[
 s_D=\frac{1}{4t(1-t)},\qquad
 z_{\mathrm{CW}}=\frac{1}{1-s_D}
 =\frac{4t(t-1)}{(2t-1)^2}.
\]
The map from the endpoint orbifold of signature $(5,5,\infty)$
to the Hecke orbifold of signature $(2,5,\infty)$ is a genuine
index-two orbifold covering.  Up to the algebraic scalar gauges,
the two RM projective period maps of the endpoint are the
restrictions of the two Cohen--Wolfart Hecke period maps
\cite[\S6.1]{CohenWolfart}.  The discrete real coordinate
uniformizes the endpoint, while the conjugate coordinate gives
its nondiscrete companion.  This statement concerns complex
projective periods and does not identify a chosen integral,
principally polarized Cohen--Wolfart moduli model.
\end{proposition}
\begin{proof}
The points $t=0,1$ map to the single order-five point, each
unramified; $t=1/2$ maps with ramification index two to the
order-two point, and the cusp at $t=\infty$ has degree two.
The orbifold Euler characteristics are respectively $-3/5$
and $-3/10$.  The index-two subgroup is the kernel of
$C_2*C_5\to C_2$ used above.
For the period comparison set $q=t(1-t)$, extend the function
field by $e^4=-q$, and put
\[
 T=\frac14\left(2-\frac{2t-1}{e^2}\right).
\]
Then $4T(1-T)=s_D$.  The substitution
$X=e^2U$, $Y=e^5V$ sends the BCDF curve $C'_5(t)$
to $V^2=U^5-5U^3+5U+2-4T$, the companion TTV curve.
The additional square root of $b(1-t)$ removes the
constant-and-variable twist of Theorem~\ref{thm:kraus-main}.
The auxiliary order-five rotations agree under this scaling,
with the RM embedding exchanged as in that theorem.
The Hecke Hilbert-modular realization of this TTV curve
also appears in \cite[\S5.2.1, (33)]{BogoLi}; its parameter
$\eta=2T-1$ gives their Hauptmodul $J=-1/s_D$.
For local solutions $f_a(s)$ of the Hecke Gauss equations
with parameters $(a,\frac12-a;1)$, the required gauges are
\[
 w_- = \left(\frac{t}{1-t}\right)^{1/4}f_{3/20}(s_D),\qquad
 w_+ = \left(\frac{t}{1-t}\right)^{-1/4}f_{1/20}(s_D).
\]
They produce the two endpoint rank-two operators, with the
discrete and conjugate embeddings assigned as above.
These identities identify the two projective period maps, one for each
real RM embedding, under the stated parameter changes and scalar gauges.
They do not specify an integral homology lattice or polarization on a
separately chosen Cohen--Wolfart factor.
\end{proof}

\begin{remark}[The Cohen--Wolfart surface and its polarization]
\label{rem:cw-polarization}
The Hecke cyclic curve
$X_z:y^{20}=x^7(1-x)^{17}(1-zx)^3$ has a primitive
abelian factor of dimension $8$, not a genus-two Jacobian;
the corrected fourth-power RM isogeny decomposition is part of
prior work \cite[\S2]{CohenWustholz}.  On the finite base
$z=c^{20}$, $T=(1-c^{-10})/2$, a degree-four curve morphism
$X_z\to C_T$ factors through two ramified double covers.
The ramification of both degree-two maps makes their Jacobian
pullbacks injective. The composite embeds $J(C_T)$ in the primitive
factor, and its image on integral first homology is saturated.
The norm identity $f_*f^*=[4]$ and adjunction show that
the ambient Jacobian polarization restricts as
$4\lambda_{J(C_T)}$, of type $(4,4)$.  Dividing this
restricted polarization by four gives the TTV principal form.
This comparison transports through the BCDF--TTV \emph{curve
isomorphism} after the displayed base extension, but the map
from the Darmon endpoint to the twisted BCDF Jacobian remains
the isogeny $\Psi$.  Corollary~\ref{cor:reverse20} identifies that endpoint with the
twisted BCDF Jacobian by $V'$, while the specified map
$\Psi=2V'$ scales its normalized principal form by $4$.
This construction does not compare the resulting homology lattice with
a separately specified Cohen--Wolfart Pochhammer-cycle lattice or
identify an associated moduli problem.  The explicit CW morphism and its
ramification are given in \cite{TTVcompanion}.
The CW comparison is restricted to $t\ne1/2$, where
$z_{\mathrm{CW}}$ has a pole.
\end{remark}

\subsection{From rational monodromy to geometric lattices}\label{sec:lattices}

The displayed matrices are integral over $R$, but their integrality
alone does not identify the homology lattice of the family. For the two
minus systems the $R$-algebra generated by monodromy determines a
unique homothety class. For the plus system there are two stable local
classes at the prime above $5$; the integral action of the quotient
involution distinguishes them. We make this comparison before reducing
modulo a coefficient prime.

\subsubsection{Stable lattice classes for the $K$-models}

We now identify the rational $K$-models of Section~\ref{sec:monodromy} with the geometric integral homology lattices.  Here \emph{saturated} has its usual integral meaning: the quotient by the indicated submodule is torsion-free.  Darmon's integral correspondences give an action of
\[
 R=\mathcal O_K=\Z[\alpha]
\]
on the torsion-free homology of the relevant abelian factors, and the connected old factors are preserved by this action \cite[Prop.~1.12 and its proof]{Darmon2000}.  In the minus cases the quotient by the connected old abelian subvariety therefore gives a saturated homology quotient.  Thus the geometric object to be identified is genuinely an $R$-lattice in the $K$-representation, not merely a $\Z$-lattice on which the displayed monodromy matrices happen to be integral.

For a given constituent write
\[
 \mathscr A=R[A,P]\subset M_2(K).
\]
A direct calculation gives
\[
 \mathscr A=R[A]+hM_2(R).
\]
Indeed $P-I=hE_{12}$ and
\[
 A(hE_{12})=hE_{22},\qquad
 (hE_{12})A=hE_{11}+xhE_{12},\qquad
 hE_{22}A=hE_{21}+xhE_{22}.
\]
Hence $hM_2(R)\subset R[A,P]$, while the reverse containment is clear from the right-hand side being a ring containing $A$ and $P$.  At every prime $\mathfrak p\nmid h$ one has
\[
 \mathscr A_{\mathfrak p}=M_2(R_{\mathfrak p}),
\]
and a lattice stable under the matrix units is of the form $I\oplus I$ for a fractional local ideal $I$.  Thus its local homothety class is unique.

For $J^-_{5,5}(1)$ one has $h=1$, so the local class is unique everywhere.  For the endpoint $J^-_{5,5}(3)$ one has $h=-2\alpha$, with $\alpha$ a unit, and the only additional prime is $(2)$.  Here
\[
 K[A]=L=\Q(\zeta_5),\qquad R[A]=\mathcal O_L.
\]
The extension $L/K$ is unramified quadratic at $(2)$, so every $A$-stable local lattice is a fractional ideal of the DVR $\mathcal O_{L,(2)}$.  Its uniformizer may be taken to be $2\in K$, and therefore all such ideals are $K_{(2)}^\times$-homothetic.  Hence the two minus systems have a unique local class at every prime.

For later use we also record the global step.  If an $R$-lattice $\Lambda\subset K^2$ is locally $a_{\mathfrak p}R_{\mathfrak p}^2$ at every prime, then $\Lambda=IR^2$ for the fractional ideal $I$ having those valuations.  Since $K=\Q(\sqrt5)$ has discriminant $5$, the Minkowski bound is $\sqrt5/2<2$; every ideal class therefore contains an integral ideal of norm $1$, so $R$ has class number one.  Consequently $I=aR$ and $\Lambda=aR^2$.

The plus system is subtler because
\[
 h=\sqrt5.
\]
Away from the unique ramified prime
\[
 \mathfrak p_5=(\sqrt5)
\]
the stable class is unique.  At $\mathfrak p_5$ one again has $R[A]=\mathcal O_L$: an eigenvalue of $A$ is a primitive fifth root, and the companion basis identifies the standard lattice with the local order $\mathcal O_{L,\mathfrak P}$.  Since $L/K$ is ramified quadratic at $5$, every $A$-stable local lattice is a fractional $\mathcal O_{L,\mathfrak P}$-ideal.  Multiplication by $K_{\mathfrak p_5}^\times$ changes the $\mathfrak P$-valuation by an even integer.  There are therefore at most two $K_{\mathfrak p_5}^\times$-homothety classes.  Both are $P$-stable, represented globally by
\[
 \Lambda_0=Re_1\oplus Re_2,
 \qquad
 \Lambda_1=Re_1\oplus R\frac{e_1-e_2}{\sqrt5}.
\]
These classes are distinct because
\[
 \Lambda_1=(A-I)^{-1}\Lambda_0,
 \qquad
 \det(A-I)=\sqrt5\,\phi,
\]
so their relative determinant has odd $\mathfrak p_5$-valuation, whereas multiplication by a scalar from $K_{\mathfrak p_5}^\times$ changes the determinant valuation by an even integer.  Equivalently, with
\[
 Q=\begin{pmatrix}1&1/\sqrt5\\0&-1/\sqrt5\end{pmatrix},
\]
the columns of $Q$ form a basis of $\Lambda_1$ and
\[
 Q^{-1}AQ=\begin{pmatrix}1&\phi\\-\sqrt5&-\phi-1\end{pmatrix},\qquad
 Q^{-1}PQ=\begin{pmatrix}1&-1\\0&1\end{pmatrix}.
\]
Monodromy stability alone therefore leaves a genuine two-class ambiguity at $5$.  We now resolve it from the integral involution.

\subsubsection{Integral $\tau$-descent for the plus family}

Work at $t_*=1/2$, so $u=-1$.  The plus equation is
\[
 y^5=u^2x^{-1}\left(\frac{x-1}{x-u}\right)^3.
\]
The point $(x,y)=(i,-1)$ is fixed by $\tau$.  To make the based topology reproducible, use
\[
 z=\frac{x-i}{x+i}.
\]
The base point becomes $z=0$, the four punctures $(-1,0,1,\infty)$ become $(i,-1,-i,1)$ in cyclic order, and $x\mapsto-1/x$ becomes $z\mapsto-z$.  Choose radial access paths and small positive meridians $(a,b,c,d)$ around these four points.  Then $abcd=1$ and the based involution is
\[
 a\mapsto c,
 \qquad
 b\mapsto d=(abc)^{-1},
 \qquad
 c\mapsto a,
 \qquad
 d\mapsto b.
\]
It squares to the identity and sends the cyclic exponents below to their negatives modulo $5$.  The chosen fixed point upstairs removes any additional deck translation.  Changing the access paths transports the entire chain model by an integral equivariant basis change; no claim that every such change is inner is needed.

Put $q=\zeta_5$, $\mathcal O_L=\Z[q]$, and use the character
\[
 (k_a,k_b,k_c,k_d)=(2,4,3,1),
 \qquad
 u_i=\frac{q^{k_i}-1}{q-1}.
\]
Each $u_i$ is a cyclotomic unit.  We briefly justify the integral compact model rather than obtaining it by scalar extension.  Let
\[
 \mathcal A=\Z[q]/(q^5-1),\qquad N=1+q+q^2+q^3+q^4.
\]
For the lifted bouquet edges $e_a,e_b,e_c$, the boundary is
\[
 \partial(e_i)=q^{k_i}-1=(q-1)u_i.
\]
Capping the finite branch points kills $Ne_i$; the cap at infinity gives the corresponding redundant $N$-linear combination.  Since
\[
 \operatorname{ann}_{\mathcal A}(q-1)=\Z N,
\]
reduction modulo $(N)$ identifies the compact homology exactly with
\[
 H=\ker\left[(u_a,u_b,u_c):\mathcal O_L^3\longrightarrow\mathcal O_L\right].
\]
There is no hidden saturation quotient.  An integral $\mathcal O_L$-basis is
\[
 f_1=\frac{e_a}{u_a}-\frac{e_b}{u_b},
 \qquad
 f_2=\frac{e_b}{u_b}-\frac{e_c}{u_c};
\]
indeed, multiplying coordinates by the units $u_i$ reduces these vectors to consecutive differences in the kernel of summation.  In the compact Gassner basis $\epsilon_i$,
\[
 f_i=(1-q)\epsilon_i.
\]

Fox differentiation of the based involution first gives the semilinear chain matrix
\[
 D_\tau=
 \begin{pmatrix}
 0&-q&1\\
 0&-q^3&0\\
 1&-q^2&0
 \end{pmatrix},
 \qquad v\longmapsto D_\tau\overline v,
\]
where the middle column comes from $b\mapsto(abc)^{-1}$.  The cap relations are preserved, and in the basis $(f_1,f_2)$ this descends to
\[
 \tau(v)=M_\tau\overline v,
 \qquad
 M_\tau=
 \begin{pmatrix}
 -1-q&-q^2-q^3\\
 -1&-q-q^2-q^3
 \end{pmatrix},
 \qquad
 M_\tau\overline{M_\tau}=I.
\]
Solving the fixed equations integrally gives
\[
 g_1=\binom q1,
 \qquad
 g_2=\binom{1+2q+2q^2+q^3}{q}.
\]
They form a full $R$-basis of $H^\tau$.  More explicitly, the integer fixed system has the four coefficients of the second coordinate as free variables, the first coordinate being an integral linear function of them; the second coordinates of
\[
 g_1,\quad \alpha g_1,\quad g_2,\quad \alpha g_2
\]
are $1,\alpha,q,\alpha q$, which form an integral basis of $\mathcal O_L$.  Thus the fixed module has the stated integral basis.

Let
\[
 p:X^+\longrightarrow C^+=X^+/\tau
\]
be the branched double cover and let $p^!$ denote homological transfer.  The two standard integral identities are
\[
 p_*p^!=2,
 \qquad
 p^!p_*=1+\tau.
\]
The first makes $p^!$ injective because $H_1(C^+,\Z)$ is torsion-free, and for every invariant cycle $v$ the second gives
\[
 2v=p^!p_*v.
\]
Hence the cokernel of transfer inside $H_1(X^+,\Z)^\tau$ is killed by $2$.  After localization at the prime above $5$, transfer is therefore an isomorphism onto the invariant lattice.  This localized statement is exactly what is needed to decide between $\Lambda_0$ and $\Lambda_1$; no global integral equality between transfer and invariant homology is used below.

\subsubsection{Transport to the $K$-monodromy coordinates}

For the character $m=1$, let $P_0=[e,T_0e]$ be the cyclotomic-to-$K$ conjugator from Section~\ref{sec:monodromy}.  On the rational fixed module put
\[
 \mathcal I=P_0^{-1}(1-q).
\]
Direct simplification gives the common nonzero scalar
\[
 \lambda=-1-2q-q^2-q^3\in L^\times
\]
with
\[
 \mathcal I(g_1)=\lambda\binom1\alpha,
 \qquad
 \mathcal I(g_2)=\lambda\binom\alpha1.
\]
Thus $\mathcal I(H^\tau\otimes\Q)=\lambda K^2$; we do not assert that $\lambda\in K^\times$.  Normalize the comparison by
\[
 J=\lambda^{-1}\mathcal I.
\]
Then $J$ is the $K$-linear comparison on the rational fixed space and
\[
 J(g_1)=\binom1\alpha,
 \qquad
 J(g_2)=\binom\alpha1.
\]
The transition determinant is
\[
 \det\begin{pmatrix}1&\alpha\\\alpha&1\end{pmatrix}
 =1-\alpha^2=\alpha\in R^\times.
\]
Hence the fixed lattice gives the class $\Lambda_0=R^2$ at $\mathfrak p_5$, not the index-$5$ class $\Lambda_1$.

The homothety ambiguity is exactly by $K^\times$.  Indeed, a matrix commuting with the nontrivial unipotent $P$ is upper triangular with equal diagonal entries, and imposing commutation with $A$ kills the remaining off-diagonal entry.  Thus the common centralizer of $A$ and $P$ in $M_2(K)$ is $K\cdot I$.  Any two $K$-linear intertwining comparisons therefore differ by a $K$-scalar, even though the cyclotomic normalization above used the scalar $\lambda\in L$.

At every prime other than $\mathfrak p_5$, the full local matrix order gives the unique lattice class; at $\mathfrak p_5$ the transfer calculation has selected $\Lambda_0$.  Applying the local-to-global ideal argument above and the class-number-one property of $R$ gives the global result.

\begin{theorem}\label{thm:lattice-final}
For each of the three $r=5$ Darmon systems, the saturated integral
homology lattice identified above, expressed in the $K$-basis of
Section~\ref{sec:monodromy}, is, up to $K^\times$-homothety,
\[
 R^2.
\]
For $J^+_{5,5}(1)$, after choosing the normalized $K$-linear comparison $J$ above, a convenient $R$-basis of the representative $R^2$ is
\[
 \binom1\alpha,\qquad \binom\alpha1.
\]
\end{theorem}

\subsubsection{The plus pairing and reduction at coefficient primes}

We now discuss the polarization only for the plus Jacobian.  The quotient $C^+$ has genus two and $\Jac(C^+)$ carries its canonical principal intersection form.  The RM action is self-adjoint for this form: on the covering Jacobian the adjoint of a deck operator is its inverse, while quotient pullback and norm are adjoint; restricting the symmetric deck correspondence to the invariant rational subspace gives the corresponding self-adjoint action on the quotient.  Under transfer, the intersection pairing is scaled by the degree $2$, so this adjointness is unaffected.

After identifying the geometric lattice with $R^2$ up to homothety, any $K$-self-adjoint rational alternating form has the shape
\[
 E(v,w)=\Tr_{K/\Q}\left(c\det_K(v,w)\right).
\]
The trace-dual ideal calculation shows that perfection on the actual geometric lattice is equivalent to
\[
 cR=\mathfrak D_K^{-1}=\frac1{\sqrt5}R.
\]
Hence
\[
 c=\frac{\varepsilon}{\sqrt5},\qquad \varepsilon\in R^\times.
\]
The unit $\varepsilon$ depends on the chosen homothety and is not determined here in the unchanged coordinates.  The choice $\varepsilon=1$ is only an illustrative representative of the inverse-different scale: on the $R$-basis $(1,\alpha)$ its cross Gram matrix is
\[
 \left(\Tr_{K/\Q}\frac{ab}{\sqrt5}\right)_{a,b=1,\alpha}
 =\begin{pmatrix}0&1\\1&-1\end{pmatrix},
\]
whose determinant is $-1$; the resulting full alternating Gram determinant is $1$.  This pairing calculation concerns only the plus Jacobian.
The principal polarization transported by $V'$ in
Corollary~\ref{cor:reverse20} and its following remark does not
determine the intersection pairings of the minus families in the
displayed $K$-bases.

Finally,
\[
 11R=(\alpha-3)(\alpha+4).
\]
Both prime ideals have norm $11$.  The lattice identification makes reduction at either prime well defined up to local integral conjugacy.  Section~\ref{sec:level11} carries out those reductions and determines the resulting projective images and three-point covers.

\begin{remark}
The conclusions above separate two phenomena that are easy to conflate. The quadratic character $\chi_{1-t}$ changes the linear lift of the minus monodromy but leaves its projective representation and stable lattices unchanged. The integral $\tau$-descent, by contrast, is essential for identifying the actual geometric lattice of the plus family at the ramified prime $5$.
\end{remark}

The integral homology lattices of the three families are isomorphic
as $R$-modules, and all six projective representations have local
orders $(5,5,\infty)$. Their real images nevertheless differ: only the
$\sigma_-$ constituent of $J^-_{5,5}(3)$ is discrete. The ordered trace
data already distinguish the cases. For the two elliptic generators,
the squared traces are equal in the endpoint minus family and unequal
in the other minus family; the rotation numbers also distinguish the
plus family. Thus the ordered monodromy data contain information
that the unordered local orders do not.

\section{Level-11 covers and their markings}\label{sec:level11}

The integral homology lattices give level structures at the coefficient
primes above $11$. Darmon's work classifies the six $r=5$ Frey systems
and their odd/even organization \cite{Darmon2000}; generation in
$\SL_2(11)$ and commutator-trace invariants also have independent
antecedents \cite{FairbairnMagaardParker,McCulloughWanderley}. We
compute the representations on these geometric lattices and the Borel
covers parametrizing their torsion lines. Although the six covers have
the same group and passport, their ordered monodromy triples differ.
For the $\sigma_-$ constituent of $J^-_{5,5}(3)$, matching the chosen
positively oriented peripheral loops, trace field, quaternionic order,
and reduction map identifies the two covers with the Clark--Voight
levels. The other four covers are distinguished from this pair in
Proposition~\ref{prop:CV11-endpoint} and
Corollary~\ref{cor:CV11-nonendpoint}.

\subsection{The two primes above 11 and the full residual image}\label{sec:residual11}

Recall
\[
 R=\Z[\alpha],\qquad \alpha^2+\alpha-1=0,
\]
so
\[
 11R=(\alpha-3)(\alpha+4).
\]
Let
\[
 \mathfrak p_1=(\alpha-3),\qquad
 \mathfrak p_2=(\alpha+4).
\]
Both residue fields are $\F_{11}$, and
\[
 \alpha\equiv3\pmod{\mathfrak p_1},\qquad
 \alpha\equiv7\pmod{\mathfrak p_2}.
\]
Reducing the geometric matrices of Section~\ref{sec:monodromy} therefore amounts simply to substituting these two values of $\alpha$. We retain the notation
\[
 A=\begin{pmatrix}0&-1\\1&x\end{pmatrix},\qquad
 B=\begin{pmatrix}x&1-xh\\-1&h\end{pmatrix},\qquad
 P=\begin{pmatrix}1&h\\0&1\end{pmatrix}.
\]
The six residue parameters and the most useful trace invariant are as follows.

\begin{center}
\begin{tabular}{c c c c c c}
\toprule
system & prime & $x$ & $h$ & $(\Tr A,\Tr B,\Tr P)$ & $\Tr[A,B]$\\
\midrule
$J^-_{5,5}(1)$ & $\mathfrak p_1$ & $3$ & $1$ & $(3,4,2)$ & $3$\\
$J^-_{5,5}(1)$ & $\mathfrak p_2$ & $7$ & $1$ & $(7,8,2)$ & $3$\\
$J^-_{5,5}(3)$ & $\mathfrak p_1$ & $3$ & $5$ & $(3,8,2)$ & $5$\\
$J^-_{5,5}(3)$ & $\mathfrak p_2$ & $7$ & $8$ & $(7,4,2)$ & $0$\\
$J^+_{5,5}(1)$ & $\mathfrak p_1$ & $7$ & $7$ & $(7,3,2)$ & $7$\\
$J^+_{5,5}(1)$ & $\mathfrak p_2$ & $3$ & $4$ & $(3,7,2)$ & $7$\\
\bottomrule
\end{tabular}
\end{center}

\begin{proposition}\label{prop:level11-image}
For each of the six rows, the subgroup generated by $A$ and $B$ is the full group
\[
 \SL_2(\F_{11}).
\]
Its projective image is $\PSL_2(\F_{11})$ of order $660$. In the minus rows the linear local orders of $(A,B,P)$ are $(5,10,11)$; in the plus rows they are $(5,5,11)$. In every case the projective orders are
\[
 (5,5,11).
\]
\end{proposition}

\begin{proof}
All identities are over $\F_{11}$. The matrices have determinant one and satisfy $ABP=I$. Since $h\ne0$ in every row, the powers of
\[
 P=U(h)=\begin{pmatrix}1&h\\0&1\end{pmatrix}
\]
run through all upper elementary matrices $U(s)$. Conjugating by $A$ gives
\[
 A U(s)A^{-1}=\begin{pmatrix}1&0\\-s&1\end{pmatrix},
\]
so the generated group contains all upper and lower elementary matrices. These generate $\SL_2(\F_{11})$. Direct enumeration gives $1320$ linear matrices and $660$ scalar classes, providing an independent check. The local orders follow from exact powers; in the minus rows $B^5=-I$, whereas in the plus rows $B^5=I$.
\end{proof}

\begin{proposition}[Generic arithmetic monodromy at $11$]
\label{prop:arithmetic-image}
Let $\mathcal A/U_K$ be any of the three abelian schemes, and
let $\mathfrak p\mid11$.  For the covariant Tate-module
convention the determinant of the corresponding arithmetic
rank-two representation is the mod-$11$ cyclotomic character:
\[
 \det\bar\rho_{\mathfrak p}=\bar\chi_{11},\qquad
 \im\bar\rho_{\mathfrak p}=\GL_2(\F_{11}),\qquad
 \im\mathbf P\bar\rho_{\mathfrak p}=\PGL_2(\F_{11}).
\]
The geometric images remain $\SL_2(\F_{11})$ and
$\PSL_2(\F_{11})$.  The resulting quadratic arithmetic
quotient is $\bar\chi_{11}^{5}$ and cuts out
$K'=K(\sqrt{-11})$.
\end{proposition}
\begin{proof}
The two real rank-two constituents are irreducible.  Each image
contains nontrivial upper and lower unipotents, and is therefore
Zariski dense in $\SL_2$.  The commutant of the full geometric
monodromy on rational first homology is $K$: the two constituents
are inequivalent because their elliptic local traces differ.
Every endomorphism defined over $K(t)$ extends over the smooth
base and commutes with this monodromy.  Since RM already supplies
$K$, it follows that $\operatorname{End}^0_{K(t)}(\mathcal A_\eta)=K$.
This assertion concerns endomorphisms over $K(t)$; no assertion
about all geometric endomorphisms is needed.  A polarization
over $K(t)$ has a Rosati involution preserving this field.  Its restriction is
the identity, since the nontrivial automorphism of a real
quadratic field is not a positive involution.
Consequently the two idempotents in
$K\otimes\Q_{11}\simeq\Q_{11}\times\Q_{11}$ are self-adjoint
for the rational Weil pairing.  The corresponding rank-two
summands are orthogonal and nondegenerate; on a two-dimensional
alternating space the similitude multiplier equals the
determinant.  The Weil pairing therefore gives
$\det\rho_{\mathfrak p,11}=\chi_{11}$, and reduction gives
the asserted determinant.  No perfect integral polarization
pairing on the minus quotients is assumed here.
By Proposition~\ref{prop:level11-image}, the geometric image
is already $\SL_2(\F_{11})$.  As
$K\cap\Q(\zeta_{11})=\Q$, the determinant character takes
all values in $\F_{11}^{\times}$ on $G_K$; hence the
arithmetic image is $\GL_2(\F_{11})$.  Passing to projective
matrices gives $\PGL_2(\F_{11})$, with determinant-square
quotient $\bar\chi_{11}^{5}$ and constant quadratic field
$K(\sqrt{-11})$.
\end{proof}

All calculations in the residual table are consequences of the displayed matrices: they are obtained by the substitutions $\alpha=3,7$ in $\F_{11}$, followed by exact matrix arithmetic.  No floating-point approximation or database recognition is used in the finite-level arguments.

The distinction between the geometric linear lift and its projectivization remains visible after reduction. If each projective generator is lifted by its unique odd-order lift, then the product of the three lifts is $-I$ in the two minus families and $+I$ in the plus family. This is exactly the persistence of Darmon's odd/even invariant in the residual triples \cite[Def.~1.4]{Darmon2000}; we do not introduce a separate ``spin'' invariant here.

\subsection{The natural degree-12 covers}\label{sec:degree12}

The group $\PSL_2(\F_{11})$ acts faithfully and transitively on
\[
 \mathbf P^1(\F_{11})=\F_{11}\cup\{\infty\},
\]
a set of $12$ points. Applying this action to the reduced triple $(A,B,P)$ gives a product-one permutation triple
\[
 (\sigma_0,\sigma_1,\sigma_\infty)\in S_{12}^3.
\]
Because an order-five element is split over $\F_{11}$, it fixes two points and has two orbits of length five. A nontrivial unipotent element of order $11$ fixes one point and acts transitively on the other eleven. Hence every row has cycle partitions
\[
 5^2 1^2,\qquad 5^2 1^2,\qquad 11\,1.
\]

\begin{proposition}[Endpoint triangular congruence covers]
\label{prop:CV11-endpoint}
Let $Y_{\mathrm{end},\mathfrak p_i}\to\mathbf P^1$ denote the degree-$12$
Borel cover obtained from the endpoint reduction at $\mathfrak p_i$.
Let $F_{\mathrm{CV}}=\mathbf Q(\lambda_5)$, where
$\lambda_n=2\cos(2\pi/n)$, and distinguish this standard trace field
from the geometric coefficient field $K=\mathbf Q(\sqrt5)$.
Put
\[
 \mathfrak q_1=(11,\lambda_5-3),\qquad
 \mathfrak q_2=(11,\lambda_5-7).
\]
Under the field identification
\[
 \iota:F_{\mathrm{CV}}\longrightarrow K,\qquad
 \lambda_{10}\longmapsto-\alpha,\qquad
 \lambda_5\longmapsto-1-\alpha,
\]
the endpoint Borel covers satisfy
\[
 Y_{\mathrm{end},\mathfrak p_i}
 \simeq_{\mathbf C}
 X_0(5,5,\infty;\mathfrak q_{3-i}),\qquad i=1,2.
\]
These are isomorphisms of covers over the ordered base
$\mathbf P^1$ with branch values $(0,1,\infty)$ and positive
peripheral generators. The endpoint projective regular covers
are likewise the corresponding principal projective congruence
covers.
\end{proposition}

\begin{proof}
We first fix the positive marking, rather than identifying the
triangle groups only abstractly. Write $\lambda=\lambda_{10}$.
Clark--Voight's extended lifts have product $-I$ and finite
traces $\lambda$, while the cusp has trace $2$
\cite[(2.1), (2.6)--(2.7)]{ClarkVoight}.
A field-defined model of this ordered triple is
\[
 D_0=\begin{pmatrix}0&1\\-1&\lambda\end{pmatrix},\qquad
 D_1=\begin{pmatrix}-\lambda&1+2\lambda^2\\-1&2\lambda\end{pmatrix},
 \qquad
 D_\infty=\begin{pmatrix}1&2\lambda\\0&1\end{pmatrix}.
\]
Here $D_0D_1D_\infty=-I$.

For an explicit comparison with their real matrices, put
$c=\cos(\pi/5)$, $s=\sin(\pi/5)$,
$v=(1+c)/s$ and $\tau=v^2>1$.
Then $\tau+\tau^{-1}=2(1+c^2)/s^2$, as required in
\cite[(2.7)]{ClarkVoight}.
The matrix
\[
 Q=\begin{pmatrix}1+2c&cv-s\\1&v\end{pmatrix},
 \qquad \det Q=\frac{2(1+c)}s>0,
\]
conjugates the source matrices
$\left(\begin{smallmatrix}c&s\\-s&c\end{smallmatrix}\right)$ and
$\left(\begin{smallmatrix}c&\tau s\\-s/\tau&c\end{smallmatrix}\right)$
to $D_0,D_1$, respectively.
Indeed $Q$ is the composition of
$z\mapsto(vz-1)/(z+v)$ and $z\mapsto c+sz$; it sends
$i,i\tau$ to $c+is,3c+is$, respectively.
At each corresponding fixed point the elliptic derivative is
$e^{2\pi i/5}$, and the traces select the displayed lifts.
The product relation then sends the cusp lift to $D_\infty$,
which is positive translation by $2\lambda$.
Thus the conjugacy lies in $\mathrm{PSL}_2(\mathbf R)$ and requires
neither inversion nor exchange of the elliptic points.
Via $\iota$ the triple becomes exactly $(-A,B,P)$; its projective
triple is therefore $( [A],[B],[P])$ in the manuscript convention.
The standard positive real embedding corresponds to $\sigma_-$.

For this signature the trace and invariant trace fields coincide:
$E=F=F_{\mathrm{CV}}$, since $\lambda_{10}=1+\lambda_5$.
The squares subgroup also equals the extended group.
In fact $\delta_j=-(\delta_j^2)^3$ for the two finite lifts,
and the cusp follows from the product relation
\cite[(5.9)--(5.12)]{ClarkVoight}.
The invariant quaternion algebra is split in the cuspidal case
\cite[Remark~5.24]{ClarkVoight}.
Let $\Phi:\mathcal B_{\mathrm{CV}}\longrightarrow M_2(K)$
be the induced algebra identification extending $\iota$.
Transporting the generated order gives
\[
 \Phi(\mathcal O_{\mathrm{CV}})=R[-A,B,P]=R[A,P].
\]
The last equality uses $B=A^{-1}P^{-1}$,
$A^{-1}=\alpha I-A$ and $P^{-1}=2I-P$.
By the geometric order identity already proved,
\[
 R[A,P]_{\mathfrak p_i}
 =R_{\mathfrak p_i}[A]+(-2\alpha)M_2(R_{\mathfrak p_i})
 =M_2(R_{\mathfrak p_i}).
\]
This is a local statement; global maximality is unnecessary.
There is no additional local distinction from Corollary~5.17,
whose annihilating factor is a unit at $11$.

Choose this transported splitting for the Clark--Voight reduction.
It sends the ordered lifts to $(-\bar A_i,\bar B_i,\bar P_i)$.
The invariant-order normalizer description gives the same
projective map: its matrices $\delta_j^2+I$ are the unit scalar
multiples $(-\alpha)(-A)$, $(-\alpha)B$, and $2P$
\cite[(5.14)--(5.23)]{ClarkVoight}.
Consequently the projective reductions agree on each positive
generator, not merely in image group.
Any other completed local splitting differs by integral
conjugacy \cite[(5.18)]{ClarkVoight}, whose residue gives a
simultaneous fractional-linear relabelling of the twelve sheets.
The cuspidal reduction is covered by Theorem~9.1 with the
authors' correction to its unipotent case
\cite[item~(8)]{ClarkVoightErrata}.

The field identification sends $\mathfrak q_2$ to $\mathfrak p_1$
and $\mathfrak q_1$ to $\mathfrak p_2$.
Taking the upper-triangular Borel in the transported splitting
therefore gives identical preimage subgroups in the marked
orbifold group; these are the Borel levels defined immediately
before \cite[(9.10)]{ClarkVoight}.
Other splitting choices give conjugate point stabilizers and
simultaneously conjugate permutation triples.
Riemann existence and smooth compactification yield the stated
marked complex isomorphisms. Equality of the projective kernels
gives the regular-cover assertion.
\end{proof}

\begin{remark}[Prime labels and the linear lift]
\label{rem:CV11-prime-swap}
The prime swap is caused by the chosen RM/trace-field identification:
at $\mathfrak p_1$, $\alpha$ reduces to $3$ and
$\iota(\lambda_5)=-1-\alpha$ reduces to $7$; at $\mathfrak p_2$
the corresponding residues are $7$ and $3$.
The geometric prime labels have not been changed.
The signed triple $(-A,B,P)$ belongs to the extended triangle
presentation with product $-I$; the geometric Jacobian lift remains
the product-one triple $(A,B,P)$.
The proposition compares projective covers over $\C$ with their ordered
branch values. Their descent to $K$ is proved next; the proposition
itself does not compare the linear regular covers.
\end{remark}

\begin{corollary}[Arithmetic descent of the endpoint Borel covers]
\label{cor:K-endpoint-borel}
Each of the six degree-$12$ marked Borel covers has exact field
of moduli $K$ and a unique model over $K$, up to unique
isomorphism over the ordered parameter line.  In the endpoint
case, the identification in Proposition~\ref{prop:CV11-endpoint}
descends to an isomorphism over $K$ with the uniquely descended
marked Clark--Voight Borel cover, with the same exchange
$\mathfrak p_i\leftrightarrow\mathfrak q_{3-i}$.
\end{corollary}
\begin{proof}
The finite $R$-module $\mathcal A[\mathfrak p]$ over $U_K$
provides the degree-$12$ torsion-line cover
\[
 \mathbf P(\mathcal A[\mathfrak p])\longrightarrow U_K.
\]
Its compactification is a marked $K$-model.  The branch-cycle power action of an element of
$G_{\Q}$ nontrivial on $K$ exchanges the two order-five
classes: their squared traces are $9$ and $5$, and powers
$\pm2$ modulo $5$ interchange them.  These classes are not
interchanged by $\PGL_2(\F_{11})$.
As proved in Proposition~\ref{prop:marked-six}, the normalizer
of the natural $\PSL_2(\F_{11})$ action in $S_{12}$ is
$\PGL_2(\F_{11})$.  Thus no relabelling of the sheets
identifies a cover with this Galois conjugate.  Its field
of moduli therefore contains $K$.
The centralizer of $\PSL_2(\F_{11})$ in its natural action
on $\mathbf P^1(\F_{11})$ is trivial: a commuting permutation
fixes the unique point fixed by a Borel and hence every point.
Thus marked cover isomorphisms are unique and their Galois
conjugates satisfy Weil descent.  The covers are algebraic over $\overline\Q$, and the
complex isomorphism of Proposition~\ref{prop:CV11-endpoint}
is defined over $\overline\Q$.  The torsion-line model
therefore supplies a $K$-descent of the marked Clark--Voight
cover.  Its uniqueness identifies it with every $K$-model
of that same marked geometric cover.  This is a comparison
in the category of covers, rather than an identification
with a separately specified abelian moduli model.
\end{proof}

\begin{remark}[Regular covers and fields of constants]
\label{rem:regular-constants}
For a $K$-defined Borel cover the arithmetic Galois closure
has group $\PGL_2(\F_{11})$, degree $1320$, and constant
field $K'=K(\sqrt{-11})$.  Over $K'$ its geometric regular
$\PSL_2(\F_{11})$ cover has degree $660$, with all labelled
deck transformations defined.  A degree-$660$ regular cover, considered without a chosen identification
of its deck group, is called a \emph{mere cover}. Such a cover can
descend to $K$: one may take the fixed field of the order-two
class of $\diag(-1,1)$ in $\PGL_2(\F_{11})$.
A regular cover with its specified
$\PSL_2(\F_{11})$-action requires $K'$.  These are fields of definition for covers and specified group actions
\cite{ClarkVoight}. They do not determine the field of moduli of the
underlying curves after the covering maps are forgotten, or give a
comparison with a separately specified polarized abelian moduli problem.
\end{remark}

\begin{corollary}[The non-endpoint marked levels]
\label{cor:CV11-nonendpoint}
The four degree-$12$ levels of $J^-_{5,5}(1)$ and $J^+_{5,5}(1)$
at the two primes above $11$ are not the two Clark--Voight covers
identified in Proposition~\ref{prop:CV11-endpoint} under
coefficient conjugation, simultaneous projective conjugacy,
peripheral inversion or exchange of the two order-five points.
\end{corollary}

\begin{proof}
The two elliptic trace-squares of each Clark--Voight reduction in Proposition~\ref{prop:CV11-endpoint} are equal.
The existing residual trace table gives $(9,5)$ and $(5,9)$ for
the generic minus system, and $(5,9)$ and $(9,5)$ for the plus
system. Equality of the two trace-squares is preserved by all
the indicated changes, including the based correction when
inversion is used as a branch-triple operation.
Thus these four marked reductions cannot be the two reductions in
Proposition~\ref{prop:CV11-endpoint}. No statement about other triangle-group
realizations, classical noncongruence, or the underlying curves
after forgetting their covering maps is implied.
\end{proof}

\begin{corollary}\label{thm:level11-covers}
Each of the six residual systems determines a connected degree-$12$ three-point cover with monodromy group $\PSL_2(\F_{11})$, passport
\[
 \boxed{(5^2 1^2;\,5^2 1^2;\,11\,1)},
\]
and genus $2$.
\end{corollary}

\begin{proof}
Transitivity follows from the full projective image. The cycle partitions above give ramification contributions
\[
 12-4=8,\qquad 12-4=8,\qquad 12-2=10.
\]
Riemann--Hurwitz therefore gives
\[
 2g-2=-2\cdot12+(8+8+10)=2,
\]
so $g=2$. Equivalently, the cover is the quotient of the regular $\PSL_2(\F_{11})$ cover by a Borel subgroup.
\end{proof}

This corollary concerns three-point covering data; it does not say that the six underlying genus-two curves are pairwise nonisomorphic when their maps are forgotten.  Their fields of moduli as marked covers are determined separately in Corollary~\ref{cor:K-endpoint-borel}; the fields of moduli of the unmarked curves need not coincide.

\begin{corollary}[Classical noncongruence and level $110$]
\label{cor:classical-noncongruence}
Under the classical identification
$\mathbf P^1\setminus\{0,1,\infty\}\simeq
\bar\Gamma(2)\backslash\mathfrak H$, the six projective
regular groups and their Borel preimages are noncongruence
subgroups of $\PSL_2(\Z)$.  All have generalized
(Wohlfahrt) level $110$.
\end{corollary}
\begin{proof}
Each base cusp of $\bar\Gamma(2)$ has width $2$.
The degree-$12$ Borel passport above yields cusp widths
\[
 (10,10,2,2),\qquad(10,10,2,2),\qquad(22,2).
\]
For the regular projective cover the cusp widths are $10$
and $22$, with multiplicities determined by the regular
inertia orbits.  In both cases their least common multiple
is $110$, but no cusp has width $110$.  If a subgroup were
congruence, the Fricke--Wohlfahrt theorem would identify its
congruence level with $110$, whereas Larcher's cusp-amplitude
criterion requires a cusp of that width
 \cite[Theorem~A and the following clarification for
$D=\Z$]{MasonSchweizer}.  Apply the criterion to the full
inverse image in $\SL_2(\Z)$, which contains $-I$; its
cusp amplitudes are precisely the displayed widths.
This contradiction proves the claim.
\end{proof}
\begin{remark}
The endpoint levels are congruence covers relative to the
nonarithmetic triangle group $\Delta(5,5,\infty)$, although
they are noncongruence subgroups of the classical modular
group.  These are different notions of congruence.
\end{remark}

\subsection{Borel levels and triangular modular curves}\label{sec:borel-triangular}

The Borel quotient parametrizes torsion lines in the specified family.
Over $U$, the integral homology local system $\Lambda$, viewed as a
rank-two $R$-module, determines
\[
 \mathbf P(\Lambda/\mathfrak p\Lambda)\simeq \mathbf P^1(\F_{11}).
\]
A point of this projective fibre is a line in the $\mathfrak p$-torsion local system, or equivalently a cyclic order-$11$ subgroup in the corresponding RM torsion component after the usual analytic comparison.  Thus the degree-$12$ cover is the level curve on which such a torsion line is chosen.  This is a level structure relative to the given abelian family;
it does not assert a universal moduli interpretation for all
hypergeometric abelian varieties.

For a \emph{discrete} projective image there is a second construction available.  Clark--Voight attach congruence covers to a triangle group by reducing its quaternionic order, and their Borel quotient is the natural triangular analogue of an $X_0$-level curve \cite{ClarkVoight}; Duque-Rosero--Voight study many low-genus examples of these covers \cite{DuqueVoight}.  The unique discrete constituent in the present paper is the $\sigma_-$ realization of $J^-_{5,5}(3)$, with image
\[
 \Delta(5,5,\infty)<\Delta(2,5,\infty)
\]
of index two.  Proposition~\ref{prop:CV11-endpoint} identifies their marked complex covers, and Corollary~\ref{cor:K-endpoint-borel} descends this agreement to the unique marked $K$-models.  Under the trace-field identification
\[
 \lambda_{10}\mapsto-\alpha,
 \qquad
 \lambda_5\mapsto-1-\alpha,
\]
the two primes are exchanged:
\[
 \mathfrak p_1\longleftrightarrow\mathfrak q_2,
 \qquad
 \mathfrak p_2\longleftrightarrow\mathfrak q_1.
\]
Thus
\[
 Y_{\mathrm{end},\mathfrak p_i}
 \simeq_{\C,(0,1,\infty)}
 X_0(5,5,\infty;\mathfrak q_{3-i}),
\]
and the projective regular covers agree over $\C$ as well.
Corollary~\ref{cor:K-endpoint-borel} strengthens the Borel comparison
to the unique marked $K$-models.  The isogeny $\Psi$ identifies the arithmetic torsion-line covers at the
exchanged RM primes: $2$ is invertible at $11$, and a scalar quadratic
twist acts trivially on lines. A comparison of polarized abelian moduli
problems would require additional compatibility of lattices and
polarizations. The torsion-line construction also does not select a
particular mere regular cover over $K$.

The four covers of $J^-_{5,5}(1)$ and $J^+_{5,5}(1)$ have unequal
elliptic trace-squares, whereas each Clark--Voight reduction in
Proposition~\ref{prop:CV11-endpoint} has equal elliptic trace-squares.
Corollary~\ref{cor:CV11-nonendpoint} therefore excludes the indicated
identifications under coefficient conjugation, simultaneous projective
conjugacy, peripheral inversion, and exchange of the order-five points.
This does not exclude every possible triangle-group realization of
these covers.

\subsection{Inequivalence of the marked covers}\label{sec:marked11}

The common passport does not identify the six reduced monodromy triples. To compare them we keep the branch values $0,1,\infty$ ordered and ask whether two triples are simultaneously conjugate. At the matrix level we consider both $\PSL_2(\F_{11})$ and $\PGL_2(\F_{11})$; at the covering level we allow an arbitrary relabeling of the twelve sheets, i.e. simultaneous conjugacy in $S_{12}$.

\begin{proposition}\label{prop:marked-six}
The six reduced monodromy triples with their branch-point orderings fixed are pairwise inequivalent under simultaneous conjugacy in
\[
 \PSL_2(\F_{11}),\qquad \PGL_2(\F_{11}),\qquad S_{12}.
\]
At either fixed prime the three Darmon systems therefore remain distinct, and for each fixed system the two coefficient primes give distinct original marked covers.
\end{proposition}

\begin{proof}
The distinction can be made conceptually from projective trace data.  With $z=\Tr(AB)=2$, the Fricke identity gives
\[
 \Tr[A,B]=\Tr(A)^2+\Tr(B)^2+z^2-
 \Tr(A)\Tr(B)z-2=2+h^2.
\]
For the six rows, the ordered invariant
\[
 \left(\Tr(A)^2,\Tr(B)^2,\Tr[A,B]\right)\in\F_{11}^3
\]
is respectively
\[
 (9,5,3),\ (5,9,3),\ (9,9,5),\ (5,5,0),\ (5,9,7),\ (9,5,7).
\]
These six triples are distinct.  The squared traces are unchanged if an $\SL_2$ lift is multiplied by $-I$, and the commutator is insensitive to those central signs; hence the displayed invariant is preserved by simultaneous conjugacy in $\PGL_2(\F_{11})$.  The six marked triples are therefore pairwise inequivalent in $\PGL_2(\F_{11})$, and a fortiori in $\PSL_2(\F_{11})$.

It remains only to pass from projective matrix conjugacy to arbitrary relabeling of the twelve sheets.  Let
\[
 G=\PSL_2(\F_{11})\leq S_{12}
\]
be its natural action on $\mathbf P^1(\F_{11})$.  If a permutation simultaneously conjugated two of our triples, then, since each triple generates $G$, that permutation would normalize $G$.  The centralizer of $G$ in this transitive coset action is trivial because a Borel subgroup is self-normalizing.  The standard automorphism theorem for $\PSL_2(q)$ therefore gives
\[
 N_{S_{12}}(G)=\PGL_2(\F_{11}):
\]
indeed $\operatorname{Out}(G)$ has order $2$ for the prime field $\F_{11}$, and the index-two extension $\PGL_2(\F_{11})$ already lies in the normalizer.  Thus an $S_{12}$ conjugator would give a $\PGL_2(\F_{11})$ conjugator, which the invariant above excludes.
\end{proof}

The commutator trace also gives a stronger distinction in the endpoint family. This is the standard Higman/commutator invariant for Nielsen equivalence of generating pairs; see \cite{McCulloughWanderley}. An automorphism of the free group on two generators sends the commutator to a conjugate of the commutator or its inverse, so its trace is unchanged. The values $5$ and $0$ attached here to the two geometric endpoint reductions therefore separate them even after changing the free generating pair and conjugating the target. The invariant itself is not new; the point is its value on these two normalized Darmon reductions.

\begin{remark}[Comparison with Chen's nonabelian level structures]
\label{rem:chen-distinction}
The six reduced epimorphisms
$F_2\twoheadrightarrow\PSL_2(\F_{11})$ also determine
points in the Hurwitz/Nielsen problem for nonabelian level
structures on punctured elliptic curves in the sense of
William Y.~Chen \cite[Table~5]{WilliamChen}.
Under inner target equivalence, their generating pairs lie in
four $\SL_2(\Z)$-components with indices
$80,72,64,80$, all of coarse genus zero.  Here inner equivalence means simultaneous conjugation
of the target epimorphism, and the mapping-class action is
by $\operatorname{Out}^+(F_2)\simeq\SL_2(\Z)$.
The central element $-I$ is absent from these stabilizers,
so the coarse degrees are $40,36,32,40$.
Chen's mapping-class stabilizers and our Borel preimages
are different subgroups; the resulting coarse genus-zero
curves are not the six genus-two torsion-line covers.  The component-size classification
belongs to the existing literature; the observation here
concerns the placement of our particular normalized pairs.
\end{remark}

\subsection{Changing the marking and the Kraus convention}

The preceding proposition keeps the branch marking fixed. There are nevertheless natural relations once the two order-five branch points are exchanged. Let
\[
 H(A,B,P)=(B,B^{-1}AB,P)
\]
be the positive Hurwitz exchange. Then the two reductions of $J^-_{5,5}(1)$ become conjugate inside $\PSL_2(\F_{11})$. For $J^+_{5,5}(1)$ the corresponding two triples become conjugate in $\PGL_2(\F_{11})$ but not in $\PSL_2(\F_{11})$; the required conjugating matrices have nonsquare determinant. By contrast, the two $J^-_{5,5}(3)$ triples remain separated by their commutator traces.

Orientation reversal should also be stated with the product convention kept explicit. Simultaneous inversion of $(A,B,P)$ need not preserve the relation $ABP=I$. A product-one version is
\[
 \mathcal R(A,B,P)=(A^{-1},B^{-1},BA).
\]
Each row is $\PGL_2(\F_{11})$-conjugate, but not $\PSL_2(\F_{11})$-conjugate, to its own $\mathcal R$-transform. Combining $H$ with $\mathcal R$ turns the two plus rows into an inner-conjugate pair. These statements describe changes of marking; they do not contradict Proposition~\ref{prop:marked-six}.

Finally, the nontrivial automorphism of $K$ exchanges the two prime substitutions. For the geometric linear local system of the $j=3$ family, the Kraus--BCDF comparison exchanges the RM embedding and introduces the character of $1-t$. Over $K(t)$ the full twist is $b(1-t)$; its constant factor has trivial complex monodromy. Thus the Kraus reduction at $\mathfrak p_1$ corresponds projectively to the Darmon $J^-_{5,5}(3)$ reduction at $\mathfrak p_2$, and conversely. The twist changes the linear triple to
\[
 (A,-B,-P),
\]
with linear orders $(5,5,22)$, but central signs disappear in the fractional-linear action, so the degree-$12$ permutation triple is unchanged.

\subsection{Relation to the companion family and further questions}

For the genus-two Jacobian in the companion study \cite{TTVcompanion},
the two real period systems have local projective orders
$(5,\infty,\infty)$. One has image $\Delta(2,5,\infty)$;
its map $s=4T(1-T)$ is a ramified pullback of projective
structures, not a faithful index-two orbifold inclusion.
The $\sigma_-$ constituent of $J^-_{5,5}(3)$ has image
$\Delta(5,5,\infty)$,
which embeds as an index-two orbifold subgroup of the Hecke
group by Proposition~\ref{prop:endpoint-hecke-periods}.
The two sets of level covers are therefore not identified
by their common Hecke source.  The algebraic curve comparison
in that proposition relates their rational periods after
finite base change and the specified quadratic twist, while
Corollary~\ref{cor:reverse20} supplies an integral comparison
between the endpoint and the BCDF family away from $2$.
The isomorphism $V'$ transports the BCDF Jacobian's principal
polarization to $J^-_{5,5}(3)$. This polarization is distinct from the
natural quotient polarization described after
Corollary~\ref{cor:reverse20}. Neither comparison specifies a
polarization on a separately chosen Cohen--Wolfart factor or identifies
the original branch markings of the companion's covers.

The arithmetic statements at level $11$ now have a precise
scope.  Proposition~\ref{prop:arithmetic-image} determines
the generic arithmetic projective image, and
Corollary~\ref{cor:K-endpoint-borel} supplies unique marked
$K$-models of the Borel covers.  The larger field required
for labelled regular covers is described in
Remark~\ref{rem:regular-constants}.  In particular, the torsion-line covers of $J^-_{5,5}(3)$ agree over
$K$, with their ordered branch values, with the descended Clark--Voight
Borel covers. A comparison with an abelian moduli problem would also
have to specify and identify its integral lattices and polarizations.

For a generalized Fermat specialization
$x^5+y^5=Cz^p$ with $x y(x^5+y^5)\ne0$, the endpoint
parameter is $t_0=x^5/S$, where $S=x^5+y^5$.
Writing $k=(xy)^2/S$, the polynomial identity
$F_{C'}(-kU,t_0)=(-k)^5F_{x,y}(U)$ gives
\[
 C'_5(t_0)\simeq C_5(x,y)^{[-k]},\qquad
 F_{x,y}(U)=U^5+5xyU^3+5x^2y^2U+y^5-x^5,
\]
in agreement with \cite[Lemma~3.5]{BCDF}.
The full endpoint twist has square class
$-k\,b(1-t_0)=-by$ in $K^\times/K^{\times2}$.
The isomorphism $V'$ therefore specializes to the Kraus
Jacobian with that twist and the same RM-prime exchange.
A rational point on a Borel torsion cover requires a $G_K$-stable line
in the chosen mod-$11$ RM constituent.  Such a line is not
forced by a solution, and generic geometric or arithmetic
surjectivity does not prove surjectivity for every fibre.
The covers are therefore not asserted to yield a new
Diophantine restriction.

Three questions remain.  The first is to determine
explicit equations and invariants of the marked genus-two
Borel covers while retaining their arithmetic descent data.
Another is to compare these $K$-models with a specified abelian moduli
problem, including its integral lattice and polarization. The embedded
TTV Jacobian in the Cohen--Wolfart factor supplies a concrete homology
lattice and restricted polarization after base extension; their
compatibility with such a moduli problem remains to be established.  A third is to determine
which higher-dimensional Frey families admit similarly
explicit integral monodromy and level comparisons, without
identifying the nondiscrete constituents with Fuchsian
triangle congruence covers.  Finite-field hypergeometric
trace formulas and the present geometric mod-$\mathfrak p$
representations concern different realizations; their full
comparison must retain both the arithmetic twists and the
local good-reduction conventions \cite{DFLST,FLRST,HoffmanTu}.

\end{document}